\documentclass[a4paper,10pt,reqno]{amsart}
\usepackage{amscd}

\usepackage{euscript}
\usepackage{epsfig}

\begin{document}

\theoremstyle{plain} \numberwithin{equation}{section}

\newtheorem{thm}{Theorem}[section]
\newtheorem{prop}[thm]{Proposition}
\newtheorem{cor}[thm]{Corollary}
\newtheorem{lemma}[thm]{Lemma}
\newtheorem{lemdef}[thm]{Lemma--Definition}
\newtheorem{ques}[thm]{Question}
\newtheorem{claim}[thm]{Claim}

\theoremstyle{definition}

\newtheorem{remark}[thm]{Remark}
\newtheorem{defn}[thm]{Definition}
\newtheorem{ex}[thm]{Example}
\newtheorem{fact}[thm]{Fact}
\newtheorem{notn}[thm]{Notation} 
\newtheorem{note}[thm]{Note}
\newtheorem{obs}[thm]{Observation}
\newtheorem{rmk}[thm]{Remark}
\newtheorem{assm}[thm]{Assumption}
\newcommand{\bi}{\begin{itemize}}
\newcommand{\ei}{\end{itemize}}
\newcommand{\bp}{\begin{proof}}
\newcommand{\ep}{\end{proof}}

\def\AA{\mathbb{A}}
\def\CC{\mathbb{C}}
\def\GGG{\mathbb{G}}
\def\PP{\mathbb{P}}
\def\QQ{\mathbb{Q}}
\def\RR{\mathbb{R}}
\def\VVV{\mathbb{V}}
\def\ZZ{\mathbb{Z}}

\def\ov{\overline}

\def\al{\alpha}
\def\be{\beta}
\def\de{\delta}
\def\eps{\epsilon}
\def\ga{\gamma}
\def\io{\iota}
\def\ka{\kappa}
\def\la{\lambda}
\def\na{\nabla}
\def\om{\omega}
\def\si{\sigma}
\def\th{\theta}
\def\ups{\upsilon}
\def\ve{\varepsilon}
\def\vp{\varpi}
\def\vt{\vartheta}
\def\ze{\zeta}

\def\De{\Delta}
\def\Ga{\Gamma}
\def\La{\Lambda}
\def\Om{\Omega}
\def\Si{\Sigma}
\def\Th{\Theta}

\def\A{\mathcal{A}}
\def\B{\mathcal{B}}
\def\C{\mathcal{C}}
\def\D{\overline{D}}
\def\Lb{\mathcal{L}}
\def\E{\overline{E}}
\def\F{\mathcal{F}}

\def\FF{\tilde{\F}}
\def\G{\mathcal{G}}
\def\GG{\tilde{\G}}
\def\HH{\mathcal{H}}
\def\I{\mathcal{I}}
\def\J{\mathcal{J}}
\def\JJ{\tilde{\J}}
\def\K{\mathcal{K}}
\def\M{\overline{M}}
\def\N{\mathcal{N}}
\def\O{\mathcal{O}}
\def\pp{\overline{p}}
\def\Q{\mathcal{Q}}
\def\S{\mathcal{S}}
\def\s{\overline{s}}
\def\S{\mathcal{S}}
\def\T{\mathcal{T}}
\def\U{\mathcal{U}}
\def\V{\mathcal{V}}
\def\VV{\tilde{\V}}
\def\W{\mathcal{W}}
\def\X{\mathcal{X}}
\def\Y{\mathcal{Y}}
\def\Z{\mathcal{Z}}

\def\Cox{\text{\rm Cox}}
\def\PGL{\text{\rm PGL}}
\def\mult{\text{\rm mult}}

\def\hdeg{\text{hdeg}}
\def\dim{\text{\rm dim}}
\def\codim{\text{codim}}
\def\Eff{\text{\rm Eff}}
\def\h{\text{h}}
\def\HH{\text{H}}
\def\H{\overline{H}}
\def\K{\text{K}}
\def\Hom{\text{Hom}}
\def\l{\overline{l}}
\def\N{\text{N}}
\def\NE{\overline{\text{NE}}}
\def\Pic{\text{\rm Pic}}
\def\Proj{\text{\rm Proj}}
\def\Spec{\text{Spec}}
\def\s{\overline{s}}
\def\Bl{\text{\rm Bl}}

\def\dra{\dashrightarrow}
\def\hra{\hookrightarrow}
\def\lra{\leftrightarrow}
\def\ra{\rightarrow}

\def\eset{\emptyset}
\def\i{\infty}
\def\setm{\setminus}

\title{The Cox Ring of $\M_{0,6}$}

\author{Ana-Maria Castravet}
\address{Department of Mathematics, University of Massachusetts at
Amherst, Amherst 01003}
\email{noni@math.umass.edu}
\date{\today}

\begin{abstract}
We prove that the Cox ring of the moduli space $\M_{0,6}$, of stable
rational curves with $6$ marked points, is finitely generated by
sections corresponding to the boundary divisors and divisors which
are pull-backs of the hyperelliptic locus in $\M_3$ via morphisms
$\rho:\M_{0,6}\ra\M_3$ that send a $6$-pointed rational curve to a
curve with $3$ nodes by identifying $3$ pairs of points. In
particular this gives a self-contained proof of Hassett and
Tschinkel's result about the effective cone of $\M_{0,6}$ being
generated by the above mentioned divisors.
\end{abstract}

\maketitle

\section{Introduction}

A question of Fulton about the moduli space $\M_{0,n}$, of stable,
$n$-pointed, rational curves, is whether the cone $\NE^k(\M_{0,n})$
of effective cycles of codimension $k$ in $\M_{0,n}$ is generated by
\emph{$k$-strata}, i.e., loci in $\M_{0,n}$ corresponding to
reducible curves with at least $k$ nodes. While the case when
$k=n-4$ (i.e., the cone of effective curves) is completely open (and
an affirmative result would imply, by results of Gibney, Keel and
Morrison \cite{GKM}, the similar statement for the moduli space
$\M_{g,n}$, of stable, $n$-pointed, genus $g$ curves, thus
determining the ample cone of $\M_{g,n}$) the case when $k=1$ (i.e.,
the cone of effective divisors) was settled independently by Keel
(unpublished, a refference to this may be found in \cite{GKM}, p.277) 
and Vermeire \cite{V}: Fulton's question has a
negative answer when $n=6$ (and therefore for any $n\geq6$). Hassett
and Tschinkel prove in \cite{HT} that the \emph{Keel-Vermeire
divisors} (pull-backs of the locus of hyperelliptic curves in the
moduli space $\M_3$, via morphisms $\M_{0,6}\ra\M_3$ sending a
$6$-pointed rational curve to a curve with $3$ nodes by identifying
$3$ pairs of points) together with the $2$-strata (the boundary)
generate the cone of effective divisors in $\M_{0,6}$. The proof in
\cite{HT} is based on a computer check. In this paper we give a
proof of Hassett and Tschinkel's result, by proving a stronger
statement: we show that the sections corresponding to the above
divisors generate the Cox ring of $\M_{0,6}$.

\

Recall that if $X$ is a smooth projective variety with Picard group
freely generated by divisors $D_1,\ldots,D_r$, then the Cox ring (or
total coordinate ring) of $X$ is the multi-graded ring:
$$\Cox(X)=\bigoplus_{(m_1,\ldots,m_r)\in\ZZ^r}\HH^0(X,m_1D_1+\ldots+m_rD_r).$$

The Cox ring being finitely generated has strong implications for
the birational geometry of $X$ ($X$ is a so-called \emph{Mori Dream Space}):
the effective cone and the nef cone are both polyhedral and there
are finitely many \emph{small modifications} of $X$ (i.e., varieties
$X'$ isomorphic in codimension one to $X$) such that any moving
divisor on $X$ (i.e., a divisor whose base locus has codimension at
least $2$) is nef on one of the varieties $X'$ (see \cite{HK} for
the precise statements). It has been conjectured by Hu and Keel
\cite{HK} that any log-Fano variety has a finitely generated Cox
ring. This has been recently proved in the groundbreaking
paper \cite{BCHM}. In \cite{HK} Hu and Keel ask the following question:

\begin{ques}\label{main}
Is the Cox ring of $\M_{0,n}$ finitely generated?
\end{ques}

As pointed out in \cite{KM}, the moduli space $\M_{0,n}$ is log-Fano
only for $n\leq6$.

\

We answer Question \ref{main} for $n=6$ by finding explicit
generators.  Our hope is that our method for finding generators,
that proved to be useful in other circumstances (see \cite{CT}),
will eventually help answer Question \ref{main} for larger $n$ as well.

\

Consider the Kapranov description of the moduli space $\M=\M_{0,6}$.
If $p_1,\ldots,p_5$ are points in linearly general position in
$\PP^3$, then $\M$ is the iterated blow-up of $\PP^3$ along
$p_1,\ldots,p_5$ and along the proper transforms of the lines
$l_{ij}=\overline{p_ip_j}$ for all $i\neq j$. If $p$ is a general
point in $\PP^3$, there is a unique twisted cubic $C$ in $\PP^3$
that contains the points $p_1,\ldots,p_5,p$. Then
$(C,p_1,\ldots,p_5,p)$  is a $6$-pointed rational curve, hence an
element of $\M_{0,6}$. The point $p$ corresponds to the $6$'th
marking (the so-called \emph{moving point}).

\

Denote by $H$ the hyperplane class on $\M$ and by $E_i$ and $E_{ij}$
the exceptional divisors in $\M$ corresponding to the points $p_i$
and the lines $l_{ij}$.

\begin{notn} Let $\La_{ijk}$ be the class of the proper
transform of the plane $\overline{p_ip_jp_k}$:
$$\La_{ijk}=H-E_i-E_j-E_k-E_{ij}-E_{ik}-E_{jk}.$$
\end{notn}

If $S\subset\{1,\ldots,6\}$ and $|S|=2,\hbox{ or }3$, let $\De_S$ be
the boundary divisor in $\M$ with general element a curve with two
irreducible components with the partition of the markings given by
$S\cup S^c$. In the Kapranov description, the boundary divisors
$\De_S$ have the following classes:
$$\De_{i6}=E_i,\quad\De_{ij6}=E_{ij}, \quad i,j=1,\ldots,5,$$
$$\De_{ij}=\La_{abc}, \hbox{ if } \{i,j,a,b,c\}=\{1,\ldots,5\}.$$

\begin{notn} Let $Q_{(ij)(kl)}$ be the class of the proper
transform of the unique quadric that contains all the points
$p_1,\ldots,p_5$ and the lines $l_{ik},l_{il},l_{jk},l_{jl}$:
$$Q_{(ij)(kl)}=2H-\sum_i E_i-E_{ik}-E_{il}-E_{jk}-E_{jl}.$$
\end{notn}

The divisor classes $Q_{(ij)(kl)}$ are exactly the divisors
considered by Keel and Vermeire: for example, if one considers the
map $\M_{0,6}\ra\M_3$ given by identifying the pairs of points
$(12)(34)(56)$, then the class of the pull-back of the hyperelliptic
locus in $\M_3$ is computed in \cite{HT} to be the class of
$Q_{(12)(34)}$. We call the divisors $Q_{(ij)(kl)}$ the
\emph{Keel-Vermeire divisors}. We prove the following:

\begin{thm}\label{main thm}
The Cox ring of  $\M_{0,6}$ is generated by the sections (unique up
to scaling) corresponding to the boundary divisors (i.e.,
$\La_{ijk}$ and the exceptional divisors $E_i$ and $E_{ij}$) and the
Keel-Vermeire divisors $Q_{(ij)(kl)}$.
\end{thm}

The paper is divided as follows: Section \ref{plan} explains the
strategy of proof, there are two main cases, the details of each are
given in Section \ref{proof main claim I}, respectively Section
\ref{proof main claim II}. The remaining sections contain auxiliary
results needed in the proof. Section \ref{inequalities} contains
proofs for some basic inequalities, while Section
\ref{multiplicities} contains some general multiplicity estimates
needed for Case II. Section \ref{proof Cox of Y} contains the proof
of Lemma \ref{Cox of Y} (needed for Case II) that states that the
Cox ring of the blow-up of $\PP^2$ in seven (non-general) points is
generated by sections corresponding to $-1$ and $-2$ curves. Section
\ref{LM} gives necessary and sufficient conditions for a divisor on
$X$, the iterated blow-up of $\PP^3$ in four general points and
lines through them, to have sections. Finally, in Section \ref{restr
to generators} we compute the restrictions of an arbitrary divisor
$D$ to all the boundary divisors and Keel-Vermeire divisors on $\M$.
Moreover, we derive some necessary conditions for these restrictions
to be effective (an assumption in our main proof).

\noindent{\em Acknowledgements.} I thank Jenia Tevelev and Sean Keel
for useful comments.

\section{Plan of Proof}\label{plan}

Consider an arbitrary divisor class on $\M$:
$$D=dH-\sum_i m_iE_i-\sum_{i,j}m_{ij}E_{ij}.$$

In all that follows we assume $\HH^0(\M,D)\neq0$.
\begin{notn}
Let $l$ be the class of the proper transform in $\M$ of a general
line in $\PP^3$. Let $e_i$ be the class of a general line in $E_i$. Let
$C$ be the class of the proper transform of a general cubic that
passes through $p_1,\ldots,p_5$:
$$C=3l-e_1-\ldots-e_5.$$
\end{notn}
The curves with class $C$ cover a dense set of $\M$; hence,
$D.C\geq0$ for any effective divisor $D$.
\begin{defn}
Let $x_{i},x_{ij},x_{ijk},x_{(ij)(kl)}$ be the sections (unique up
to scalar) corresponding  to the divisors:
\begin{equation}\label{generators}
E_i,E_{ij},\La_{ijk},Q_{(ij)(kl)}.
\end{equation}
\end{defn}

\begin{defn}
We call a section $s\in\HH^0(\M,D)$ a \emph{distinguished section}
if
$$s=x_i^{n_i}x_{ij}^{n_{ij}}x_{ijk}^{n_{ijk}}x_{(ij)(kl)}^{n_{(ij)(kl)}},$$
where $n_i,n_{ij},n_{ijk},n_{(ij)(kl)}$ are non-negative integers.
\end{defn}

To show that $\HH^0(\M, D)$ is generated by distinguished sections,
we do an induction on $D.C$. Note that we may assume that $D$
contains none of the divisors (\ref{generators}) in its base locus,
i.e., equivalently, if for $E$ any of the divisors in
(\ref{generators}), one has $\HH^0(E,D_{|E})\neq0$. To see this,
note that if $E$ is an effective divisor, say $E$ is the zero locus
of a section $x_E\in\HH^0(\M,E)$, then there is an exact sequence:
$$0\ra\HH^0(\M, D-E)\ra\HH^0(\M, D)\ra\HH^0(E, D_{|E}).$$

If $\HH^0(E, D_{|E})=0$ then any $s\in\HH^0(\M,D)$ is of the form
$x_Et$, where $t\in\HH^0(\M,D-E)$. If in addition $E$ is a divisor
in (\ref{generators}) then we may replace $D$ with $D-E$ and $s$
with $t$. (Clearly, if $t$ is generated by distinguished sections,
then $s$ is too.) Therefore, we may assume:

\begin{assm}\label{assm}
$\HH^0(E,D_{|E})\neq0$ for all divisors $E$ in (\ref{generators}).
\end{assm}

Denote by $r_E$ the restriction to $E$:
\begin{equation*}
r_E:\HH^0(\M,D)\ra\HH^0(E,D_{|E}).
\end{equation*}

To prove Theorem \ref{main thm} it is enough to prove the following:

\

\noindent {\bf Main Claim.} Let $D$ be a divisor on $\M$:
$$D=dH-\sum_i m_iE_i-\sum_{i,j}m_{ij}E_{ij},$$
such that $\HH^0(\M,D)\neq0$ and that satisfies Assumption
\ref{assm}. Up to a renumbering, we may assume that $m_5\leq m_i$,
for $i=1,\ldots,4$. If $m_i=m_5$ for all $i=1,\ldots,4$, then we may
assume that the maximum of the $m_{ij}$'s for all
$i,j\in\{1,\ldots,5\}$ is attained for $m_{i5}$ for some
$i=1,\ldots,4$. Let $E=E_5$. Then for any $s\in\HH^0(D\M,D)$, there
is $s'\in\HH^0(\M,D)$, generated by distinguished sections, such
that $r_E(s)=r_E(s')$.

\

To see how the Main Claim implies Theorem \ref{main thm}, note that
the kernel of the restriction $r_E$ is $\HH^0(\M,D-E)$ and the map
$\HH^0(\M,D-E)\ra\HH^0(\M,D)$ is given by multiplication with $x_E$.
If $r_E(s)=r_E(s')$, then $s-s'=x_Et$, where $t\in\HH^0(\M,D-E)$. If
$s'$ is generated by distinguished sections, then to show that $s$
is generated by distinguished sections it is enough to show that
$\HH^0(\M,D-E)$ is generated by distinguished sections. We may
replace $D$ with $D-E$, and continue the procedure. Since $E$ is
always among the $E_i$'s, note that $(D-E).C<D.C$ and
$\HH^0(\M,D-E)$ is generated by distinguished sections by induction.
The process has to stop as $D.C\geq0$ for any effective divisor $D$.
(In particular, note that $D.C$ decreases also when we substract
from $D$ any of the divisors $E$ in (\ref{generators}) for which
$\HH^0(E,D_{|E})=0$.)

\begin{notn}\label{D bar}
Given any divisor $D$ on $\M$ we denote by $\D$ the restriction
$D_{|E_5}$ of $D$ to $E_5$. By (\ref{restr to E_i}) one has:
\begin{equation*}
\D=m_5\H-\sum_{i=1}^4m_{i5}\E_i.
\end{equation*}
\end{notn}

Let $\rho_5:\PP^3\dra\PP^2$ be the projection from $p_5$. Let
$q_i=\rho_5(p_i)$ ($i\in\{1,\ldots, 4\}$). The divisor $E_5$ is
isomorphic to the blow-up of $\PP^2$ along the points
$q_1,\ldots,q_4$ (as $q_i$ determines the direction of the line
$l_{i5}$). The divisors $\H$, respectively $\E_i$, are the
hyperplane class, respectively the exceptional divisors on $E_5$
(see also Section \ref{section restr to E_i}.) The map $\rho_5$ is
resolved by the morphism $\pi_5:\M\ra\M_{0,5}$ that forgets the
$5$'th marking (which is also a retract for the inclusion
$E_5\subset\M$).

\begin{notn}
Let $\l_{ij}$ be the line $\overline{q_iq_j}$ in $\PP^2$. Denote:
$$x=\l_{13}\cap\l_{24},\quad
y=\l_{14}\cap\l_{23},\quad z=\l_{12}\cap\l_{34}.$$
\end{notn}

\begin{notn}
Let $L_x$ be the proper transform in $\M$ of the unique line in
$\PP^3$ that passes through $p_5$ and intersects the skew lines
$l_{13}$ and $l_{24}$. Similarly, let $L_y$ (respectively $L_z$) be
the unique line that passes through $p_5$ and intersects the skew
lines $l_{14}$ and $l_{23}$ (respectively $l_{12}$ and $l_{34}$).
\end{notn}
Remark that $x=\rho_5(L_x)$, $y=\rho_5(L_y)$, $z=\rho_5(L_z)$.

\

In order to prove the Main Claim, we distinguish two cases.

\

\noindent {\bf \underline{Case I:  Assume that $D.L_x\geq0,\quad
D.L_y\geq0,\quad D.L_z\geq0$.}}

\

\begin{notn}\label{sections}
Denote by  $s_{ij}$ the section on $E_5$ corresponding to the proper
transform of the line $\l_{ij}$ in $\PP^2$. Let $s_i$
($i=1,\ldots,4$) be the sections corresponding to the exceptional
divisors $\E_i$.
\end{notn}

\begin{defn}
We call a section $s\in\HH^0(E_5,\D)$ a \emph{distinguished section
on $E_5$} if $s$ can be written as a monomial in the sections
$s_{ij}$ and $s_i$.
\end{defn}

Since $E_5\cong\M_{0,5}$ is the blow-up of $\PP^2$ along
$q_1,\ldots,q_4$, by Lemma \ref{eff cone} the Cox ring $\Cox(E_5)$
of $E_5$ is generated by distinguished sections. The Main Claim
follows from the following:

\begin{prop}\label{case I}
Under the assumptions of the Main Claim and the assumptions in Case
I, the restriction map
$$r_{E_5}:\HH^0(\M,D)\ra\HH^0(E_5,\D),
$$ is surjective and one may
lift any distinguished section (hence, any section) in
$\HH^0(E_5,\D)$ to a section generated by distinguished sections in
$\HH^0(\M,D)$.
\end{prop}

The following is the main observation needed to prove Proposition
\ref{case I}:

\

\noindent {\bf Main Observation -- Case I.} Distinguished sections
on $E_5$ may be lifted to distinguished sections on $\M$ using the
following rules:
\begin{equation}\label{rules I}
{x_{ij5}}_{|E_5}=s_{ij},\quad {x_{i5}}_{|E_5}=s_i.
\end{equation}
This is because ${\La_{ij5}}_{|E_5}=\l_{ij}$, ${E_{i5}}_{|E_5}=\E_i$
(see Section \ref{restr to generators}, formula (\ref{restr to
E_i})).

\noindent {\bf Sketch of Proof of Proposition \ref{case I}.} We lift
a distinguished section $\s\in\HH^0(E_5,\D)$ using the rules
(\ref{rules I}). Hence, there is a section $t'$ belonging to some
$\HH^0(\M,D')$, where $\D'=\D$ and $r_{E_5}(t')=\s$.

\begin{notn}
Let $\De=D-D'$.
\end{notn}

% correction 1
\begin{notn}\label{X}
Denote by $X$ the iterated blow-up of $\PP^3$ in $p_1,\ldots,p_4$
and proper transforms of lines $l_{ij}$ ($i,j\in\{1,\ldots,4\}$).
\end{notn}

Since $D$ and $D'$ have the same restriction to $E_5$, it follows
from (\ref{restr to E_i}) that the divisor $\De$ is a divisor on $X$.
% end correction 1
Note, $X$ is a toric variety. The following is a standard result:

\begin{lemma}\label{toric}
The Cox ring of $X$  is generated by sections $x_i$,$x_{ij}$,
$x_{ijk}$ corresponding to the exceptional divisors $E_i$,$E_{ij}$
and proper transforms of hyperplanes $\La_{ijk}$
($i,j,k\in\{1,\ldots4\}$).
\end{lemma}

% correction 2
Proposition \ref{case I} is now  immediate if $\HH^0(\De)\neq0$:
Since the points $p_1,\ldots, p_5$ are general, the restriction to $E_5$ of any distinguished section in $\Cox(X)$ is non-zero. In particular, 
if $t''$ is any non-zero  section in $\HH^0(\De)$, then $t''_{|E_5}\in\HH^0(E_5,\O)$ is non-zero.
Therefore, the section $s=t't''$ is a section in $\HH^0(\M,D)$ that restricts to (a non-zero multiple of) $\s$ in $\HH^0(E_5,\D)$. 
Since $t''$ is a distinguished section, it follows that $t$ is
generated by distinguished sections.
% end correction 2

\begin{defn}
We call a distinguished section $\s$ on $E_5$ \emph{a section with
straightforward lifting to $D$} if after lifting using the rules
(\ref{rules I}) we end up with a divisor $D'$ for which $\De=D-D'$
has $\HH^0(\De)\neq0$.
\end{defn}

The following Claim (proof in Section \ref{proof main
claim I}) finishes the proof of Proposition \ref{case I}.

\begin{claim}\label{main claim I}
Under the assumptions of Proposition \ref{case I}, any distinguished
section $\s\in\HH^0(E_5,\D)$ is a linear combination of
distinguished sections with straightforward lifting to $D$.
\end{claim}

\

{\bf \underline{Case II:  Assume one of $D.L_x, D.L_y, D.L_z$ is
negative.}}

\

\begin{defn}\label{m's}
Let:
$$m_x=\hbox{max }\{0,-D.L_x\},\quad m_y=\hbox{max }\{0,-D.L_y\},
\quad m_z=\hbox{max }\{0,-D.L_z\}.$$
\end{defn}

\begin{notn}
Denote by  $Y$ the blow-up of $\PP^2$ along $q_1,q_2,q_3,q_4,x,y,z$.
Let $\E_i$,$\E_x$,$\E_y$,$\E_z$ be the corresponding exceptional
divisors. For a given divisor $D$ on $\M$ we consider the following
divisor $\D^Y$ on $Y$:
$$\D^Y=\D-m_x\E_x-m_y\E_y-m_z\E_z.$$
\end{notn}

Clearly, the linear system $\HH^0(Y, \D^Y)$ is a subspace of the
linear system $\HH^0(E_5,\D)$.

\begin{claim}\label{factors}
The restriction map $r_{E_5}$ factors through $\HH^0(Y, \D^Y)$.
\end{claim}

\bp Clearly, Claim \ref{factors} is non-trivial only when one of
$m_x$,$m_y$,$m_z$ is positive. Take for example the case when
$m_x>0$ (the other cases are identical).  By Proposition
\ref{multiplicity estimates}, the line $L_x$ is  contained in $D$
with multiplicity $m\geq m_x>0$. It follows that for any
$s\in\HH^0(\M,D)$ the section $r_{E_5}(s)$ vanishes at $x$ with
multiplicity $\geq m$; hence, $r_{E_5}(s)$ lies in the subspace
$\HH^0(Y, \D^Y)$. \ep

In Case II we follow the exact same steps as in Case I, with the
only difference that we work on $Y$ instead of $E_5$.

\begin{notn}
Denote by  $s'_{ij}$ the section corresponding to the proper
transform in $Y$ of the line $\l_{ij}$.  Similarly, let $s_{xy},
s_{xz}, s_{yz}$ be the sections corresponding to the proper
transforms of the lines $\overline{xy}, \overline{xz},
\overline{yz}$.  Let $s_i, s_x, s_y,s_z$ be the sections
corresponding to the exceptional divisors $E_i, E_x,E_y,E_z$.
\end{notn}

Note:
$$x\in\l_{13},\l_{24},\quad y\in\l_{14},\l_{23},\quad z\in\l_{12},\l_{34}.$$
Hence, for example $s'_{13}$ is a section of the divisor
$\H-\E_1-\E_3-\E_x$ and  the section $s_{13}$ (Notation
\ref{sections}) is given by $s_{13}=s'_{13}s_x$.  Moreover, if we
let:  $$r_Y:\HH^0(\M,D)\ra\HH^0(Y,\D^Y),$$ be the morphism of Claim
\ref{factors}, then $r_{E_5}(s)=r_Y(s)s_x^{m_x}  s_y^{m_y}
s_z^{m_z}$.

\begin{defn}
We call a section $s\in\HH^0(Y,\D^Y)$ a \emph{distinguished section
on $Y$} if $s$ can be written as a monomial in the sections
$s'_{ij}$,$s_{xy}$,$s_{xz}$,$s_{yz}$, $s_i$,$s_x$,$s_y$,$s_z$.
\end{defn}

In Section \ref{proof Cox of Y}  we prove the following:
\begin{lemma}\label{Cox of Y}
The Cox ring $\Cox(Y)$ of $Y$ is generated by distinguished
sections.
\end{lemma}

Note, by Lemma \ref{Cox of Y}, the generators of $\Cox(Y)$ are given
by the sections (unique up to scalar multiplication) corresponding
to the $(-1)$ and $(-2)$ curves on $Y$. The Main Claim follows from:

\begin{prop}\label{case II}
Under the assumptions of the Main Claim, the
restriction map:
$$r_Y:\HH^0(\M,D)\ra\HH^0(Y,\D'),$$ is surjective and one may
lift any distinguished section (hence, any section) in
$\HH^0(Y,\D^Y)$ to a section generated by distinguished sections in
$\HH^0(\M,D)$.
\end{prop}

The following is the main observation needed to prove Proposition
\ref{case II}:

\

\noindent {\bf Main Observation -- Case II.} Distinguished sections
on $Y$ may be lifted to distinguished sections on $\M$ using the
following rules:
\begin{equation}\label{rules IIa}
r_Y(x_{ij5})=s'_{ij},\quad  r_Y(x_{i5})=s_i,
\end{equation}
\begin{equation}\label{rules IIb}
r_Y(x_{(13)(24)})=s_{yz},\quad r_Y(x_{(14)(23)})=s_{xz},\quad
r_Y(x_{(12)(34)})=s_{xy}.
\end{equation}

This is because when $D=\La_{ij5}$ one has:
$$\D^Y=\H-\E_i-\E_j-\E_{\al},$$ where $\al=x$ if $ij\in\{13,24\}$,
$\al=y$ if $ij\in\{14,23\}$, $\al=z$ if $ij\in\{12,34\}$. Similarly:
$$\overline{Q}_{(13)(24)}^Y=\H-\E_y-\E_z,
\overline{Q}_{(14)(23)}^Y=\H-\E_x-\E_z,
\overline{Q}_{(12)(34)}^Y=\H-\E_x-\E_y.$$

\noindent {\bf Sketch of Proof of Proposition \ref{case II}.} We
lift a distinguished section $\s\in\HH^0(E_5,\D)$ using the rules
(\ref{rules IIa}) and (\ref{rules IIb}). Hence, there is a section
$t'$ in some $\HH^0(\M,D')$ and $r_Y(t')=\s$. As in Case I, we let
$\De=D-D'$. The divisor $\De$ is a divisor on $X$ (Notation
\ref{X}). As in Case I, Proposition \ref{case II} follows from Lemma
\ref{toric} if $\HH^0(\De)\neq0$.

\begin{defn}
We call a distinguished section $\s$ on $Y$ \emph{a section with
straightforward lifting to $D$} if lifting using the rules
(\ref{rules IIa}) and (\ref{rules IIb}) results in a divisor $D'$
for which $\De=D-D'$ has $\HH^0(\De)\neq0$.
\end{defn}

The following Claim (proof in Section \ref{proof main
claim II}) finishes the proof of Proposition \ref{case II}.

\begin{claim}\label{main claim II}
Under the assumptions of the Main Claim, any distinguished section
$\s\in\HH^0(Y,\D^Y)$ is a linear combination of distinguished
sections with straightforward lifting to $D$.
\end{claim}

\section{Proof of Claim  \ref{main claim I}}\label{proof main claim I}

The idea is that any distinguished section on $E_5$ can be
rewritten, using the relations in $\Cox(E_5)$, as a linear
combination of distinguished sections with straightforward lifting.
To check that $\HH^0(\De)\neq0$ we use Lemma \ref{LMspaces}.
Assumption \ref{assm} is equivalent to inequalities (\ref{ineq
E_i}), (\ref{ineq E_ij}), (\ref{ineq planes}), (\ref{ineq quadrics})
(for all permutations of indices).

We use the notation from Section \ref{LM}. Recall that $e_{ij}$ is
the class of a fiber of the $\PP^1$-bundle $E_{ij}\ra\l_{ij}$. One
has $D.l=d$, $D.e_i=m_i$, $D.e_{ij}=m_{ij}$ (see for example
(\ref{basic restr to E_i}), (\ref{basic restr to E_ij})). The
inequalities defining Case I are equivalent to:
\begin{equation}\label{extra lines}
D.(l-e_5-e_{ij}-e_{kl})\geq0, \quad\{i,j,k,l\}=\{1,2,3,4\}.
\end{equation}

\begin{lemma}\label{lift I}
Let $\s$ be a distinguished section on $E_5$:
\begin{equation}\label{simple s}
\s=\prod_{i,j}s_{ij}^{a_{ij}}\prod_{i}s_i^{l_i},
\end{equation}
where $a_{ij},l_i\geq0$. If $\s$ is a section $\HH^0(E_5,\D)$ then
$\s$ has straightforward lifting to $D$ if and only if for all
$\{i,j,k,l\}=\{1,2,3,4\}$ one has:

\begin{equation}\label{ineq}
 a_{ij}\leq D.(C_{k;l}-e_5),
\end{equation}
where $C_{k;l}=2l-e_{ki}-e_{kj}-e_l$.
\end{lemma}

\begin{rmk}\label{nonneg}
By (\ref{ineq quadrics}) one has $D.(C_{k;l}-e_5)\geq0$ for all
$k,l\in\{1,2,3,4\}$.
\end{rmk}

\begin{rmk}
The condition that $\s$ is in $\HH^0(\D)$ is equivalent to:
\begin{equation}\label{sum I}
\sum a_{ij}=D.e_5,\quad a_{ij}+a_{ik}+a_{il}-l_i=D.e_{i5},
\end{equation}
(the coefficients of $\H$ and $\E_i$ in $\D$). It follows from
(\ref{sum I}) that:
\begin{equation}\label{rel}
a_{kl}-a_{ij}-l_k-l_l=D.(e_{k5}+e_{l5}-e_5),
\end{equation}
\begin{equation}\label{sum of l's}
\sum_{i=1}^4l_i=D.(2e_5-\sum_{i\neq 5}e_{i5}),
\end{equation}
\begin{equation}\label{cyclic sum}
a_{jk}+a_{jl}+a_{kl}+l_i=D.(e_5-e_{i5}).
\end{equation}
\end{rmk}

\bp[Proof of Lemma~\ref{lift I}] If $\D=0$, $\s=1$ (i.e., $a_{ij}=0,
l_i=0$) then the lift $D'$ is $0$. Hence, $\De=D-D'=D$. Since
$\HH^0(D)\neq0$, there is nothing to prove in this case.

Assume now $\D\neq0$. Recall that $E_5\subset\M$ has a retract
$\pi:\M\ra E_5\cong\M_{0,5}$ given by the morphism that forgets the
$5$-th marking. One has:
\begin{equation}\label{pull backs}
\pi^*\l_{ij}=\La_{ij5}+E_{ij},\quad\pi^*\E_i=\E_{i5}+E_i
\end{equation}
(This is a general fact about the forgetful morphisms
$\pi_i:\M_{0,n}\ra\M_{0,n-1}$ that forget a marking $i$. If $\De_S$
is a boundary divisor in $\M_{0,n-1}$, corresponding to the
partition $S\cup S^c$, then $\pi^*\De_S=\De_{S}+\De_{S\cup\{i\}}$.)

Since we lift $\D$ to $D'$ by lifting $\l_{ij}$ to $\La_{ij5}$ and
$\E_i$ to $E_{i5}$, it follows that:
$$D'=\pi^*\D-\De_0,$$
where $\De_0$ is the effective divisor on $X$ given by:
$$\De_0=\sum_{i,j\in\{1,\ldots,4\}}a_{ij}E_{ij}+
\sum_{i\in\{1,\ldots,4\}}l_i E_i.$$

Then $\De=D-D'=D-\pi^*\D+\De_0$.

\begin{obs}\label{obs1}
If $(D-\pi^*\D).C\geq0$ for some nef curve $C$ on $X$, then
$\De.C\geq0$.
\end{obs}

Below we show that $(D-\pi^*\D).C\geq0$ for all the nef curves $C$
in Lemma \ref{LMspaces} giving inequalities $(1)$--$(4)$. Hence, by
Observation \ref{obs1}, $\De.C\geq0$. For the remaining nef curves
$C$ in general it will not be true that $(D-\pi^*\D).C\geq0$, but we
show that we still have $\De.C\geq0$ for the nef curves $C$ giving
inequalities $(5),(7),(8),(9)$ and that for $C=C_{k;l}$ (inequality
$(6)$) $\De.C\geq0$ is equivalent to (\ref{ineq}). Note:
\begin{equation}
(\pi^*\D).l=D.e_5,\quad(\pi^*\D).e_i=D.e_{i5},\quad(\pi^*\D).e_{ij}=0\quad
(i,j\neq5)
\end{equation}
(It is enough to check this when $D=H,E_i,E_{ij}$. For this, use the
formulas (\ref{pull backs}).)

We check one by one the inequalities $(1)-(9)$ in Lemma \ref{LMspaces}:
$$(1)\quad (D-\pi^*\D).l=D.(l-e_5)\geq0,$$
as $l-e_5$ is a nef curve on $\M$. Similarly:
$$ (2)\quad (D-\pi^*\D).(l-e_i)=D.(l-e_i-e_5+e_{i5})\geq0
\hbox{ by }(\ref{ineq E_ij}),$$

$$(3)\quad (D-\pi^*\D).(l-e_{ij})=D.(l-e_5-e_{ij})\geq0
\hbox{ by } (\ref{ineq planes}),$$

$$(4)\quad (D-\pi^*\D).(l-e_{ij}-e_{kl})=D.(l-e_5-e_{ij}-e_{kl})\geq0
\hbox{ by } (\ref{extra lines}).$$

For inequality $(5)$ (recall $C_{ij}=2l-e_{ij}-e_k-e_l$):
$$ (D-\pi^*\D).C_{ij}=D.(C_{ij}-2e_5+e_{k5}+e_{l5}),$$
$$\De_0.C_{ij}=a_{ij}+l_k+l_l,$$
$$\De.C_{ij}=D.(C_{ij}-e_5)+D.(e_{k5}+e_{l5}-e_5)+a_{ij}+l_k+l_l.$$
By (\ref{rel}), $\De.C_{ij}=D.(C_{ij}-e_5)+a_{kl}$. From (\ref{ineq
planes}) (and $a_{kl}\geq0$), $\De.C_{ij}\geq0$.

For inequality $(7)$ (recall $C_i=2l-e_{ij}-e_{ik}-e_{il}$):
$$(D-\pi^*\D).C_i=D.(C_i-2e_5).$$
and $\De_0.C_i=a_{ij}+a_{ik}+a_{il}$. Using (\ref{sum I}),
$\De_0.C_i=D.e_{i5}+l_i$. Therefore:
$$\De.C_i=D.(C_i-2e_5+e_{i5})+l_i=2D.(l-e_i-e_5+e_{i5})+D.(2e_i-\sum_{u\neq i}e_{iu})+l_i.$$
It follows from  (\ref{ineq E_i}) and (\ref{ineq E_ij}) that $\De.C_i\geq0$.

For inequality $(8)$ (recall $B=3l-\sum_{i=1}^4e_i$):
$$(D-\pi^*\D).B=D.(B-3e_5+\sum_{i=1}^4e_{i5}),$$
and $\De_0.B=\sum_{i=1}^4l_i$.
It follows from (\ref{sum of l's}) that $\De.B=D.(3l-\sum_{i=1}^5e_i)\geq0$.

For inequality $(9)$ (recall $B_i=3l-2e_i-e_{jk}-e_{jl}-e_{kl}$):
$$(D-\pi^*\D).B_i=D.(B_i-3e_5+2e_{i5}),$$
and $\De_0.B_i=a_{jk}+a_{jl}+a_{kl}+2l_i$. From (\ref{cyclic sum})
one has $\De_0.B_i=D.(e_5-e_{i5})+l_i$,
$$\De.B_i=D.(B_i-2e_5+e_{i5})=
2D.(l-e_i-e_5+e_{i5})+D.(l-e_{jk}-e_{jl}-e_{kl}-e_{i5}).$$ It
follows by (\ref{ineq E_ij}) and (\ref{ineq planes}) that
$D.B_i\geq0$.

There is at least one strict inequality in $(4)$: assume
$\De.(l-e_{ij}-e_{kl})=0$, for all $\{i,j,k,l\}=\{1,2,3,4\}$. From
the computation above for case $(4)$ we have:
$$(D-\pi^*\D).(l-e_{ij}-e_{kl})=D.(l-e_5-e_{ij}-e_{kl})\geq0.$$
As $\De_0.(l-e_{ij}-e_{kl})=0$ ($l-e_{ij}-e_{kl})$ is a nef curve)
it follows that:
$$(D-\pi^*\D).(l-e_{ij}-e_{kl})=\De_0.(l-e_{ij}-e_{kl})=0.$$
Since $\De_0.(l-e_{ij}-e_{kl})=a_{ij}+a_{kl}$ it follows that
$a_{ij}=0$ for all $i,j$. By (\ref{sum I}), $D.e_5=0$ and
$D.e_{i5}=0$,$l_i=0$ for all $i\neq5$. Hence, $\D=0,\s=1$, which
contradicts our assumption.

We show now that inequality $(6)$ is equivalent to (\ref{ineq}). One
has:
$$(D-\pi^*\D).C_{i;j}=D.(C_{i;j}-2e_5+e_{j5}),$$
and $\De_0.C_{i;j}=a_{ik}+a_{il}+l_j$. From (\ref{cyclic sum}),
$\De_0.C_{i;j}=D.(e_5-e_{j5})-a_{kl}$. Therefore:
$$\De.C_{i;j}=D.(C_{i;j}-e_5)-a_{kl}.$$
Hence, inequality $(6)$ is equivalent to (\ref{ineq}). \ep

\subsection {\bf Proof of Claim \ref{main claim I}.} Let $\s$ be a
distinguished section in $\HH^0(E_5,\D)$ as in (\ref{simple s}). If
$a_{ij}\leq D.(C_{k;l}-e_5)$ for all $\{i,j,k,l\}=\{1,2,3,4\}$ then
by Lemma \ref{lift I} $\s$ has straightforward lifting to $D$.
Assume now that $a_{ij}>D.(C_{k;l}-e_5)$ for some choice of
$i,j,k,l$. Without loss of generality, we may assume
$a_{12}>D.(C_{3;4}-e_5)$. Note that by Remark \ref{nonneg} it
follows that $a_{12}>0$.

\begin{claim}\label{claim}\label{positivity}
If $a_{12}>D.(C_{3;4}-e_5)$ then either $a_{34}>0$ or $l_1+l_2>0$.
\end{claim}

\bp By (\ref{rel}) one has:
\begin{equation}\label{a_12-a_34}
a_{12}-a_{34}-l_1-l_2=D.(e_{15}+e_{25}-e_5).
\end{equation}

Assume $a_{34}=l_1=l_2=0$. It follows from (\ref{a_12-a_34}) and
$a_{12}> D.(C_{3;4}-e_5)$ that
$$a_{12}=D.(e_{15}+e_{25}-e_5)>D.(C_{3;4}-e_5).$$
This is a contradiction, as by (\ref{ineq quadrics}) one has:
\begin{equation*}
D.(C_{3;4}-e_5)-D.(e_{15}+e_{25}-e_5)=
 D.(2l-e_4-e_{13}-e_{23}-e_{15}-e_{25})\geq0.
\end{equation*}
\ep

\subsection{\bf Algorithm for replacing $\s$.}\label{alg I}
We now give an algorithm for replacing $\s$ with another
distinguished section $\s'$ for which $a_{12}-D.(C_{3;4}-e_5)$ is
strictly smaller than for $\s$ and moreover, for all $i,j$ for which
$a_{ij}-D.(C_{k;l}-e_5)$ increases by this change, the section $\s'$
(still) satisfies $a_{ij}-D.(C_{k;l}-e_5)\leq0$. We repeat the
following two steps until $a_{12}\leq D.(C_{3;4}-e_5)$ (as by Claim
\ref{claim} one of the two situations must happen if
$a_{12}>D.(C_{3;4}-e_5)$). The same argument works for any $a_{ij}$.

\

\underline{Step 1: If $l_1+l_2>0$}: We may assume without loss of
generality that $l_1>0$. Consider the following sections in the
linear system $|\H-\E_2|$:
$$s_{12}s_1,\quad s_{23}s_3,\quad s_{24}s_4.$$

The linear system $|\H-\E_2|$ is $1$-dimensional and any two of the
above sections are linearly independent. Since $a_{12}>0,l_1>0$, we
may replace $s_{12}s_1$ in $s$ with a linear combination of
$s_{23}s_3$ and $s_{24}s_4$. The effect on the coefficients $a_{ij}$
and $l_i$ (of the corresponding two distinguished sections) is as
follows: $a_{12}$ and $l_1$ both decrease by $1$, while either
$a_{23},l_3$ increase by $1$, or $a_{24},l_4$ increase by $1$
(everything else stays the same). But by Lemma \ref{no adjacent
bads} one has:
$$a_{2j}<D.(C_{k;l}-e_5),\quad\hbox{ for all } j\in\{3,4\},
\{j,k,l\}=\{1,3,4\}.$$ Therefore, after
increasing $a_{23}$ or $a_{24}$ by $1$ one still has $a_{2j}\leq
D.(C_{k;l}-e_5)$.

\

\underline{Step 2: If $a_{34}>0$}: Consider the following sections
in the linear system $|2\H-\E_1-\ldots-\E_4|$:
$$s_{12}s_{34}, \quad s_{13}s_{24},\quad  s_{14}s_{23}.$$

The linear system $|2\H-\E_1-\ldots-\E_4|$ is $1$-dimensional and
any two of the above sections are linearly independent. Since
$a_{12}>0,a_{34}>0$, we may replace $s_{12}s_{34}$ in $s$ with a
linear combination of $s_{13}s_{24}$ and $s_{14}s_{23}$. The effect
on the coefficients $a_{ij}$ is: $a_{12}$ and $a_{34}$ both decrease
by $1$, while either $a_{13},a_{24}$ increase by $1$, or
$a_{14},a_{23}$ increase by $1$. By Lemma \ref{no adjacent bads} one
has:
$$a_{ij}<D.(C_{k;l}-e_5),\quad\hbox{ for all }i\in\{1,2\}, j\in\{3,4\},
\{i,j,k,l\}=\{1,2,3,4\}.$$

Therefore after increasing $a_{13},a_{14},a_{23},a_{24}$ by $1$,
each of them still satisfies its corresponding inequalities.

\begin{lemma}\label{no adjacent bads}
If $a_{ij}> D.(C_{k;l}-e_5)$ then $a_{iu}<D.(C_{v;w}-e_5)$ for all
$\{u,v,w\}=\{j,k,l\}$ such that $u\in\{k,l\}$.
\end{lemma}

\bp Assume the contrary. Then $a_{ij}+a_{iu}>
D.(C_{k;l}-e_5)+D.(C_{v;w}-e_5)$. But by (\ref{cyclic sum})
$a_{ij}+a_{iu}\leq D.(e_5-e_{u'5})$, where $\{u',u\}=\{k,l\}$. This
is a contradiction with Claim \ref{sums of C_k;l}. \ep

\begin{claim}\label{sums of C_k;l}
$D.(C_{k;l}+C_{v;w}-2e_5)\geq D.(e_5-e_{u'5})$ for all $v,w,u'$ such
that $\{u',u\}=\{k,l\}$ and $\{v,w,u\}=\{j,k,l\}$ for some
$u\in\{k,l\}$.
\end{claim}

\bp There are four cases:

Case (i): $v=j,w=l$ ($u=k$,$u'=l$). Using (\ref{ineq E_ij}) and
(\ref{ineq planes}) one has:
\begin{gather*}
D.(C_{k;l}+C_{j;l}-2e_5)-D.(e_5-e_{l5}) =\\
=2D.(l-e_5-e_l+e_{l5})+D.(l-e_{ij}-e_{ik}-e_{jk}-e_{l5})+D.(l-e_5-e_{jk})\geq0.
\end{gather*}

Case (ii): $v=l,w=j$ ($u=k$,$u'=l$). Using (\ref{ineq E_ij}) and
(\ref{ineq planes}) one has:
\begin{gather*}
D.(C_{k;l}+C_{l;j}-2e_5)-D.(e_5-e_{l5})=D.(l-e_5-e_l+e_{l5})+\\
D.(l-e_5-e_j+e_{j5})+D.(l-e_{iu}-e_{il}-e_{ul}-e_{j5})+D.(l-e_5-e_{jk})\geq0.
\end{gather*}

Case (iii): $v=j,w=k$ ($u=l$,$u'=k$). This is symmetric to Case
(ii).

Case (iv): $v=k,w=j$ ($u=l$,$u'=k$). Using (\ref{ineq E_i}),
(\ref{ineq E_ij}) and (\ref{ineq planes}) one has:
\begin{gather*}
D.(C_{k;l}+C_{k;j}-2e_5)-D.(e_5-e_{l5})=\\
=2D.(l-e_5-e_k+e_{k5})+D.(2l-e_5-e_j-e_l-e_{ik})+D.(2e_k-\sum_{\al\neq
k}e_{k\al})\geq0.
\end{gather*}
 \ep

\section{Proof of Claim  \ref{main claim II}}\label{proof main claim II}

As in Section \ref{proof main claim I}, we show that any
distinguished section on $Y$ can be rewritten, using the relations
in $\Cox(Y)$, as a linear combination of distinguished sections with
straightforward lifting. Assumption \ref{assm} is equivalent to the
inequalities (\ref{ineq  E_i}), (\ref{ineq E_ij}), (\ref{ineq
planes}), (\ref{ineq quadrics}) (for all permutations of indices).
We use the  notation from Section \ref{LM}.

\begin{notn}
Let $\chi:\{12,13,14,23,24,34\}\ra\{x,y,z\}$ be the function
$$\chi(13)=\chi(24)=x,\quad\chi(14)=\chi(23)=y,\quad\chi(12)=\chi(34)=z.$$
\end{notn}

Note, one has:
$$L_{\chi(ij)}=l-e_5-e_{ij}-e_{kl},\quad\hbox{ for all }
\{i,j,k,l\}=\{1,2,3,4\}.$$

\begin{rmk}\label{D.L}
By Definition \ref{m's} one has $m_{\al}+D.L_{\al}\geq0$ for all
$\al\in\{x,y,z\}$, with equality if and only if $D.L_{\al}\leq0$.
\end{rmk}

\begin{lemma}\label{lift II}
Let $\s$ be a distinguished section on $Y$:
\begin{equation}\label{s}
\s=\prod_{i,j}{s'}_{ij}^{a_{ij}}\prod_{i}s_i^{l_i}
s_{xy}^{c_z}s_{xz}^{c_y}s_{yz}^{c_x}s_x^{l_x}s_y^{l_y}s_z^{l_z},
\end{equation}
where $a_{ij},l_i,c_x,c_y,c_z,l_x,l_y,l_z\geq0$. If $\s$ is a
section $\HH^0(Y,\D^Y)$ then $\s$ has straightforward lifting to $D$
if and only if for all $\{i,j,k,l\}=\{1,2,3,4\}$ and
$\al\in\{x,y,z\}$:
\begin{equation*}\tag{i}\label{i}
\quad c_{\chi(ij)}-a_{ij}\leq D.(l-e_5-e_{ij}),
\end{equation*}
\begin{equation*}\tag{ii}\label{ii}
c_{\al}-l_{\al}\leq  m_{\al}+D.L_{\al},
\end{equation*}
\begin{equation*}\tag{iii}\label{iii}
c_{\chi(ij)}-a_{ij}\leq D.(C_{kl}-e_5),
\end{equation*}
\begin{equation*}\tag{iv}\label{iv}
a_{ij}+\sum_{\al\neq\chi(ij)}c_{\al}\leq D.(C_{k;l}-e_5),
\end{equation*}
\begin{equation*}\tag{v}\label{v}
c_x+c_y+c_z\leq D.C,
\end{equation*}
where $C_{kl}=2l-e_i-e_j-e_{kl}$, $C_{k;l}=2l-e_{ki}-e_{kj}-e_l$,
$C=3l-\sum_{i=1}^5e_i$.
\end{lemma}

\begin{rmk}\label{RHS}
Note that the right sides of the inequalities in Lemma \ref{lift II}
are non-negative due to (\ref{ineq planes}) (for (\ref{i})),
(\ref{iii})), (\ref{ineq quadrics}) (for (\ref{iv})),  Remark
\ref{D.L} (for (\ref{ii})) and because $C$ is a nef curve on $\M$
(for (\ref{v})).
\end{rmk}

\begin{rmk}\label{comb cond}
The condition that $\s$ is in $\HH^0(Y,\D^Y)$ is equivalent to:
\begin{equation}\label{sum}
\sum a_{ij}+(c_x+c_y+c_z)=D.e_5,
\end{equation}
\begin{equation}\label{sum1}
a_{ij}+a_{ik}+a_{il}-l_i=D.e_{i5},
\end{equation}
(the coefficients of $\H$ and $\E_i$ in $\D^Y$)
\begin{equation}\label{m_x}
a_{ij}+a_{kl}+\sum_{\al\neq\chi(ij)}c_{\al}-l_{\chi(ij)}=m_{\chi(ij)},
\end{equation}
(the coefficient of $\E_\al$ in $\D^Y$ for $\al\in\{x,y,z\}$). From
(\ref{sum}), (\ref{sum1}) and (\ref{m_x}) one has:
\begin{equation}\label{sum2}
\sum l_i+2(c_x+c_y+c_z)=D.(2e_5-\sum_{i\neq5}e_{i5}),
\end{equation}
\begin{equation}\label{sum3}
a_{jk}+a_{jl}+a_{kl}+(c_x+c_y+c_z)+l_i=D.(e_5-e_{i5}),
\end{equation}
\begin{equation}\label{sum4}
(c_x+c_y+c_z)-(l_x+l_y+l_z)=(m_x+m_y+m_z)-D.e_5.
\end{equation}
\end{rmk}

\bp[Proof of Lemma~\ref{lift II}] We lift $\s$ using the rules
(\ref{rules IIa}) and (\ref{rules IIb}) (see also Remark \ref{enough
s_x}) to a section of the divisor:
\begin{gather*}
D'=\sum a_{ij}\La_{ij5}+\sum l_i E_{i5}+c_xQ_{(13)(24)}+c_yQ_{(14)(23)}+c_zQ_{(12)(34)}=\\
=(\sum  a_{ij}+2\sum c_{\al})H-\sum_{i\neq5}(a_{ij}+a_{ik}+a_{il}+\sum c_{\al})E_i-(\sum a_{ij}+\sum c_{\al})E_5-\\
-\sum_{i,j\neq5}(a_{ij}+\sum_{\al\neq\chi(ij)}
c_{\al})E_{ij}-\sum_{i\neq5}(a_{ij}+a_{ik}+a_{il}-l_i)E_{i5}.
\end{gather*}

Using (\ref{sum}) and (\ref{sum1}) one has:
\begin{gather*}
D'=(D.e_5+\sum c_{\al})H-\sum_{i\neq5}(D.e_{i5}+l_i+\sum c_{\al})E_i-(D.e_5)E_5-\\
-\sum_{i,j\neq5}(a_{ij}+\sum_{\al\neq\chi(ij)}
c_{\al})E_{ij}-\sum_{i\neq5}(D.e_{i5})E_{i5}.
\end{gather*}

Then $\De=D-D'$ is given by the following formula:
\begin{gather*}
\De=(D.(l-e_5)-\sum c_{\al})H-\sum_{i=1}^4 (D.(e_i-e_{i5})-l_i-\sum
c_{\al})E_i-\\
-\sum_{i,j\neq5}(D.e_{ij}-a_{ij}-\sum_{\al\neq\chi(ij)}
c_{\al})E_{ij}.
\end{gather*}

We show that $\De.C\geq0$ for the nef curves $C$ giving the
inequalities $(1),(2),(7),(9)$ in Lemma \ref{LMspaces} and that for
the nef curves $C$ giving the remaining inequalities, $\De.C\geq0$
is equivalent to
(\ref{i}),(\ref{ii}),(\ref{iii}),(\ref{iv}),(\ref{v}).

\

For inequality $(1)$:
$$\De.l=D.(l-e_5)-\sum c_{\al}.$$
By (\ref{sum3}), one has $\sum c_{\al}\leq D.(e_5-e_{i5})$. By the
assumption in the Main Claim $D.e_5\leq D.e_i$. Then $\De.l\geq
D.(l-e_5-e_i+e_{i5})$. It follows from  (\ref{ineq E_ij}) that
$\De.l\geq0$.

\

For inequality $(2)$:
$$\De.(l-e_i)=D.(l-e_5-e_i+e_{i5})+l_i.$$
It follows from (\ref{ineq E_ij}) that $\De.(l-e_i)\geq0$.

\

Inequality $(3)$ is equivalent to (\ref{i}) as one has:
$$\De.(l-e_{ij})=c_z-a_{12}\leq D.(l-e_5-e_{ij})+a_{ij}-c_{\chi(ij)}.$$

\

For inequality $(4)$:
$$\De.(l-e_{ij}-e_{kl})=D.(l-e_5-e_{ij}-e_{kl})+
a_{ij}+a_{kl}+\sum_{\al\neq\chi(ij)}-c_{\chi(ij)}.$$ By using
(\ref{m_x}) to substitute $a_{ij}+a_{kl}+\sum_{\al\neq\chi(ij)}$ one
has that $\De.(l-e_{ij}-e_{kl})\geq0$ is equivalent to (\ref{ii}).
Note that in Lemma \ref{LMspaces} we require that at least one of
the inequalities is strict. As Lemma \ref{one strict} shows, this is
automatically satisfied in this case.

\

For inequality $(5)$:
$$\De.C_{kl}=D.(C_{kl}-2e_5+e_{i5}+e_{j5})+a_{kl}+l_i+l_j-c_{\chi(kl)}.$$
Using (\ref{sum1}) (to substitute $l_i, l_j$) and (\ref{sum})
$\De.C_{kl}\geq0$ is equivalent to $(iii)$.

\

For inequality $(6)$:
$$\De.C_{k;l}=D.(C_{k;l}-2e_5+e_{l5})+a_{ik}+a_{jk}+l_l+c_{\chi(ij)}.$$
By using (\ref{sum3}) to substitute
$a_{ik}+a_{jk}+l_l+c_{\chi(ij)}$, $\De.C_{k;l}\geq0$ is equivalent
to $(iv)$.

\

For inequality $(7)$ (recall that $C_i=2l-e_{ij}-e_{ik}-e_{il}$):
$$\De.C_i=D.(C_i-2e_5)+a_{ij}+a_{ik}+a_{il}.$$
By using (\ref{sum1}) to substitute $a_{ij}+a_{ik}+a_{il}$,
$\De.C_i=D.(C_i-2e_5+E_{i5})+l_i$. But:
$$D.(C_i-2e_5+E_{i5})=2D.(l-e_i-e_5+m_{i5})+
D.(2e_i-\sum_{j\neq i}e_{ij}).$$
From (\ref{ineq E_ij}) and (\ref{ineq E_i}) it follows that $\De.C_i\geq0$.

\

For inequality $(8)$ (recall that $B=3l-\sum_{i=1}^4e_i$):
$$\De.B=D.(B-3e_5+\sum_{i\neq 5}e_{i5})+\sum l_i+\sum c_{\al}.$$
By using (\ref{sum2}) to substitute $\sum l_i+2\sum c_{\al}$,
$\De.B\geq0$ is equivalent to $(v)$.

\

For inequality $(9)$ (recall that $B=3l-2e_i-e_{jk}-e_{jl}-e_{kl}$):
$$\De.B_i=D.(B_i-3e_5+2e_{i5})+a_{jk}+a_{jl}+a_{kl}+2l_i+\sum c_{\al}.$$
By using (\ref{sum3}) to substitute $a_{jk}+a_{jl}+a_{kl}+l_i+\sum
c_{\al}$, $\De.B_i=D.(B_i-2e_5+e_{i5})$. But one has:
$$D.(B_i-2e_5+e_{i5})=2D.(l-e_i-e_5+e_{i5})+D.(l-e_{jk}-e_{jl}-e_{kl}-e_{i5}).$$
It follows from (\ref{ineq E_ij}) and (\ref{ineq planes}) that $\De.B_i\geq0$.
\ep

\begin{rmk}\label{enough s_x}
In order to lift $\s\in\HH^0(Y,\D^Y)$ we need to group $s'_{ij}$
with $s_{\chi(ij)}$, such that we may lift
$s_{ij}=s'_{ij}s_{\chi(ij)}$ to $x_{ij5}$, etc (so in fact we lift
$\s s_x^{m_x}s_y^{m_y}s_z^{m_z}$). For this we need to have enough
sections $s_x$,$s_y$,$s_z$. Take the case of $s_x$: one needs
exactly $a_{13}+a_{24}+c_y+c_z$ of them (to be distributed to
$s_{13}$,$s_{24}$,$s_y$,$s_z$). Since the image of the restriction
map $r_Y$ in $\HH^0(E_5,\D)$ is
$$\HH^0(Y,\D')s_x^{m_x}s_y^{m_y}s_z^{m_z},$$ the number of
$s_x$'s appearing in $\s s_x^{m_x}s_y^{m_y}s_z^{m_z}$ is  $m_x+l_x$
and by (\ref{m_x}) one has:
$$m_x+l_x=a_{13}+a_{24}+c_y+c_z.$$
\end{rmk}

\begin{lemma}\label{one strict}
It is not possible to have $c_{\al}-l_{\al}\geq  m_{\al}+D.L_{\al}$
for all $\al\in\{x,y,z\}$.
\end{lemma}

\bp Assume the contrary and add up the three inequalities. Then one
has:
$$\sum c_{\al}-\sum l_{\al}\geq \sum m_{\al}+\sum D.L_{\al}.$$
By (\ref{sum4}), this is equivalent to $\sum D.L_{\al}\leq-D.e_5$,
which contradicts Lemma \ref{sum of m_x's}. \ep

\subsection {\bf Proof of Claim \ref{main claim II}.} Let $\s$ be a
distinguished section in $\HH^0(Y,\D^Y)$ as in (\ref{s}). If
inequalities (\ref{i})-(\ref{v}) in Lemma \ref{lift II} are
satisfied, then by Lemma \ref{lift II} $\s$ has straightforward
lifting to $D$. Assume now that one of the inequalities
(\ref{i})-(\ref{v}) fails. We first show that we can keep replacing
the section $\s$ with a linear combination of distinguished sections
until we are in one of the following cases:
\begin{gather*}
(A)\quad c_x=c_y=c_z=0,\\
(B)\quad l_x=l_y=l_z=0, c_x+c_y+c_z>0,\\
(C)\quad c_x=c_y=l_x=l_y=0, c_z>0,l_z>0\quad(\hbox{up to a
permutation of $x,y,z$}).
\end{gather*}

This follows from:
\begin{claim}\label{replacements}
If $l_{\al}>0$ and $c_{\be}>0$ for $\al,\be\in\{x,y,z\}$,
$\be\neq\al$, then we may replace $\s$ with a sum of distinguished
sections $\s'$ for which both $c_x+c_y+c_z$ and $l_x+l_y+l_z$
decreased.
\end{claim}

\bp We may assume without loss of generality that $l_x>0,c_z>0$.
Consider the following sections in the linear system $|\H-\E_y|$:
$$s_{xy}s_x,\quad s'_{14}s_1s_4,\quad s'_{23}s_2s_3.$$

The linear system $|\H-\E_y|$ is $1$-dimensional and any two of the
above sections are linearly independent. Hence, we may replace
$s_{xy}s_x$ with a linear combination of the sections
$s'_{14}s_1s_4, s'_{23}s_2s_3$. The effect is: $c_z,l_x$ decrease by
$1$ and either $a_{14},l_1,l_4$ or $a_{23},l_2,l_3$ increase by $1$.
Note that $c_x,c_y,l_y,l_z$ stay the same. Hence, both $c_x+c_y+c_z$
and $l_x+l_y+l_z$ decreased by $1$.
\ep

Note that while doing replacements as in Claim \ref{replacements} we
ignore how the changes affect inequalities in Lemma \ref{lift II}.

\subsection{\underline{{\bf Case (A): $c_x=c_y=c_z=0$}}} This case is
very similar to Case I.

\begin{lemma}\label{lift IIA}
If $c_x=c_y=c_z=0$ then $\s$ has straightforward lifting to $D$
if and only if for all $\{i,j,k,l\}=\{1,2,3,4\}$ one has:
\begin{equation*}
a_{ij}\leq D.(C_{k;l}-e_5).
\end{equation*}
\end{lemma}

\bp One may immediately see (use for example Remark \ref{RHS}) that
the inequalities  (\ref{i}),(\ref{iii}),(\ref{v}) in Lemma \ref{lift
II} are satisfied. The inequality (\ref{ii})  is satisfied (see
Remark \ref{D.L}). Condition (\ref{iv}) in Lemma \ref{lift II}
becomes $a_{ij}\leq D.(C_{k;l}-e_5)$ in Case (A). \ep

\noindent {\bf Algorithm for replacing $s$ -- Case (A).} If for all
$\{i,j,k,l\}=\{1,2,3,4,\}$ one has $a_{ij}\leq D.(C_{k;l}-e_5)$ and
by Lemma \ref{lift IIA} $\s$ has straightforward lifting to $D$. if
for some $i,j,k,l$ one has $a_{ij}>D.(C_{k;l}-e_5)$ we will replace
$\s$ with a sum of distinguished sections such that all the
inequalities improve, while leaving $c_x=c_y=c_z=0$. We do this in
exactly the same way as we did in Case I, as Lemma \ref{no adjacent
bads}, Claim \ref{positivity}, as well as the Algorithm \ref{alg I}
all apply word by word.

\subsection{\underline{{\bf Case (B): $l_x=l_y=l_z=0$,
$c_x+c_y+c_z>0$}}} This is impossible because of (\ref{sum4}) and
Lemma \ref{sum of m_x's}.

\subsection{\underline{{\bf Case (C): $c_x=c_y=l_x=l_y=0$,
$c_z>0,l_z>0$}}}

\begin{rmk}
Under the assumptions of Case (C) the relations in Remark
\ref{comb cond} become:
\begin{gather}
\sum a_{ij}+c_z=D.e_5\label{sum adjusted},\\
a_{13}+a_{24}+c_z=m_x,\quad a_{14}+a_{23}+c_z=m_y\label{m_x adjusted},\\
a_{12}+a_{34}-l_z=m_z\label{m_z adjusted},\\
c_z-l_z=m_x+m_y+m_z-D.e_5\label{sum4 adjusted}.
\end{gather}
\end{rmk}

From  (\ref{m_x adjusted}) one has:
\begin{equation}\label{m_x>0}
0<c_z\leq \text{ min }\{m_x, m_y\}.
\end{equation}

From the definitions of $m_x,m_y$ it follows that
$m_x=-D.L_x$,$m_y=-D.L_y$. From (\ref{m_z adjusted}) and (\ref{sum4
adjusted}) one has:
\begin{equation}\label{a_12+a_34}
a_{12}+a_{34}-c_z=m_5-m_x-m_y=D.(2l-e_5-e_{13}-e_{14}-e_{23}-e_{24}).
\end{equation}

\begin{lemma}\label{lift IIC}
Under the assumptions of Case (C) $\s$ has straightforward lifting
to $D$ if and only if:
\begin{equation*}\tag{iii'}\label{iii'}
a_{ij}\leq D.(C_{ij}-e_5)+D.(2l-e_5-e_{13}-e_{14}-e_{23}-e_{24}),
\end{equation*}
\begin{equation*}\tag{iv'}\label{iv'}
a_{ij}\leq D.(C_{k;l}-e_5),
\end{equation*}
whenever either $ij=12, kl=34$ or $ij=34, kl=12$.
\end{lemma}

\begin{rmk}\label{RHS'}
By (\ref{ineq planes}) and (\ref{ineq quadrics}) the right hand
sides of (\ref{iii'}), (\ref{iv'}) are $\geq0$.
\end{rmk}

\bp[Proof of Lemma~\ref{lift IIC}] We claim that in Lemma \ref{lift
II} the inequalities (\ref{i}),(\ref{ii})  and (\ref{v}) are
satisfied and that (\ref{iii}), respectively (\ref{iv}) are
equivalent to (\ref{iii'}) and (\ref{iv'}).

Inequality (\ref{i}): by Remark \ref{RHS} the inequalities involving
$c_x,c_y$ are automatic. We claim that $c_z\leq D.(l-e_5-e_{ij})$
whenever $ij=12$ or $34$: by (\ref{m_x>0}) one has $c_z\leq
m_x,m_y$, hence $c_z\leq(m_x+m_y)/2$ and the claim follows from
Lemma \ref{m_x+m_y}.

Inequality (\ref{ii}): this is clearly satisfied for $l_x-c_x=0$,
$l_y-c_y=0$. From (\ref{sum4 adjusted}) and Lemma \ref{sum of m_x's}
it follows that $c_z-l_z\leq0$ and we are done by Remark \ref{RHS'}.

Inequality (\ref{iii}): the inequalities involving $c_x$ and $c_y$
are automatically satisfied. The inequalities (\ref{iii}) involving
$c_z$ are of the form (here $ij=12$ or $34$):
\begin{equation}\label{above}
c_z-a_{ij}\leq D.(C_{kl}-e_5).
\end{equation}

Using (\ref{a_12+a_34}) to substitute $c_z-a_{ij}$ in (\ref{above}),
one obtains (\ref{iii'}):
$$\quad a_{kl}\leq D.(C_{kl}-e_5)+D.(2l-e_5-e_{13}-e_{14}-e_{23}-e_{24}).$$

Inequality (\ref{iv}): We claim that the inequalities involving
$a_{13},a_{14},a_{23},a_{24}$ are satisfied: this is because by
(\ref{m_x adjusted}) $a_{ij}+c_z\leq m_x$ whenever $ij\neq12, 34$.
By Lemma \ref{ineq m_x} $m_x\geq D.(C_{k;l}-e_5)$ and we are done.
The inequalities (\ref{iv}) involving $a_{12},a_{34}$ are exactly
the inequalities (\ref{iv'}).

Inequality (\ref{v}): this follows from (\ref{m_x>0}) and Lemma
\ref{m_x+m_y}. \ep

\subsection{\bf Algorithm for replacing $s$ in Case (C)}
If the inequalities in Lemma \ref{lift IIC} are satisfied, then $\s$
has straightforward lifting to $D$. Assume one of (\ref{iii'}) or
(\ref{iv'}) is not satisfied, say for $a_{12}$ (the same argument
applies for $a_{34}$). Then by Remark \ref{RHS'} one has $a_{12}>0$.
Then we make replacements to decrease $a_{12}$ as follows: Consider
the following sections in the linear system
$|2\H-\E_1-\E_2-\E_x-\E_y|$:
$$s'_{12}s_{xy}s_z,\quad s'_{13}s_{23}s_3^2,\quad s'_{14}s_{24}s_4^2.$$

The linear system is $1$ dimensional and any two of the above
sections are linearly independent. Since $a_{12},c_z,l_z>0$, we may
replace $s'_{12}s_{xy}s_z$ in $s$ with a linear combination of
$s'_{13}s'_{23}s_3^2,s'_{14}s'_{24}s_4^2$. The effect is:
$a_{12},c_z,l_z$ decrease by $1$, while either $a_{13},a_{23}$
increase by $1$, or $a_{14},a_{24}$ increase by $1$. Note that
besides the above changes and the changes affecting the $l_i$'s
(which we ignore, since they do not appear in (\ref{iii'}),
(\ref{iv'}) no other changes occur. In particular, we still have
$c_x=c_y=l_x=l_y=0$.

\

The inequalities involving $a_{12}$ were improved (while
the ones involving $a_{34}$ remained the same). If after the
replacement $c_z=0$ or $l_z=0$, we are in Case (A) or Case (B), we
apply the procedure described for those cases. If after the
replacement we still have $c_z>0$ and $l_z>0$, then we are in Case
(C) and therefore all inequalities are satisfied, except perhaps
(\ref{iii'}), (\ref{iv'}) for $a_{12}$ or $a_{34}$.

\section{Inequalities involving $m_x$,$m_y$,$m_z$}\label{inequalities}

The assumptions in this section are the same as in the Main Claim.
Recall:
$$L_{\chi(ij)}=L_{\chi(kl)}=l-e_5-e_{ij}-e_{kl}.$$

\begin{lemma}\label{ineq m_x}
For any $\{i,j,k,l\}=\{1,2,3,4\}$ one has:
$$-D.L_{\chi(kl)}\leq D.(C_{k;l}-e_5),$$
where  $C_{k;l}=2l-e_l-e_{ik}-e_{jk}$.
\end{lemma}

\bp One has:
\begin{gather*}
D.(C_{k;l}-e_5)+D.L_{\chi(kl)}=D.(l-e_{ik}-e_{jk}-e_{ij}-e_{l5})+\\
+D.(l-e_5-e_l+e_{l5})+D.(l-e_5-e_{kl})\geq0.
\end{gather*}
It follows from (\ref{ineq E_ij}) and (\ref{ineq planes}) that
$D.(C_{k;l}-e_5)+D.L_{\chi(kl)}\geq0$. \ep

\begin{lemma}\label{m_x+m_y}
For any $i,j\in\{1,2,3,4\}$ one has:
$$-\frac{1}{2}D.(\sum_{\al\neq\chi(ij)}L_{\al})\leq\hbox{min
}\{d-m_5-m_{ij},3d-\sum_{i=1}^5m_i\}$$.
\end{lemma}

\bp Without loss of generality, we may assume $ij=12$. One has:
\begin{gather*}
2D.(l-e_5-e_{12})+D.(L_x+L_y)=\\
=2D(l-e_5-e_{12})+D.(2l-2e_5-e_{13}-e_{14}-e_{23}-e_{24})=\\
=D.(l-e_{12}-e_{13}-e_{23}-e_{45})+(l-e_{12}-e_{14}-e_{24}-e_{35})+\\
+D.(l-e_5-e_3+e_{35})+D.(l-e_5-e_4+e_{45})+D.(e_3+e_4-2e_5).
\end{gather*}
The first inequality follows from (\ref{ineq planes}) and the
assumption $D.e_5\leq D.e_i$. Moreover:
\begin{gather*}
2D(3l-\sum_{i=1}^5e_i)+D.(L_x+L_y)=\\
2(3l-\sum_{i=1}^5e_i)+D.(2l-2e_5-e_{13}-e_{14}-e_{23}-e_{24})=\\
D.(2l-e_1-e_3-e_5-e_{24})+D.(2l-e_1-e_4-e_5-e_{23})+\\
D.(2l-e_2-e_3-e_5-e_{14})+D.(2l-e_2-e_4-e_5-e_{13}).
\end{gather*}
The second inequality now follows from (\ref{ineq planes}). \ep

\begin{lemma}\label{sum of m_x's}
One has $(m_x+m_y+m_z)\leq D.e_5$,$-D.(L_x+L_y+L_z)< D.e_5$.
\end{lemma}

\bp Note, by definition of $m_x$ if $m_x>0$, then $m_x=-D.L_x$
(similarly for $y,z$). If $m_x=m_y=m_z=0$. The Claim follows from
(\ref{ineq E_i}).

Case 1) Assume just one of $m_x,m_y,m_z$ is $>0$, say $m_x>0,
m_y=m_z=0$:
\begin{equation*}
D.e_5-(m_x+m_y+m_z)=D.(l-e_{13}-e_{24}).
\end{equation*}
But $D.(l-e_{13}-e_{24})\geq0$ (see Lemma \ref{LMspaces}). The other
cases are similar.

Case 2) Assume two of $m_x,m_y,m_z$ is $>0$, say $m_x,m_y>0,m_z=0$:
\begin{equation*}
D.e_5-(m_x+m_y+m_z)=D.(2l-e_5-e_{13}-e_{14}-e_{23}-e_{24}).
\end{equation*}
By (\ref{ineq quadrics}) $D.e_5-(m_x+m_y+m_z)\geq0$. The other cases
are similar.

Case 3) Assume $m_x,m_y,m_z>0$:
\begin{gather*}
D.e_5-(m_x+m_y+m_z)=D.(L_x+L_y+L_z+e_5)=D.(3l-2e_5-\sum_{i,j=1,\ldots
4}e_{ij})=\\=D.(2e_i-\sum_{j\neq i}e_{ij})+2D(l-e_5-e_i+e_{i5})
+D.(l-e_{jk}-e_{kl}-e_{jl}-e_{i5}),
\end{gather*}
for any $\{i,j,k,l\}=\{1,2,3,4\}$. By (\ref{ineq E_i}),(\ref{ineq
E_ij}),(\ref{ineq planes}) $D.e_5-(m_x+m_y+m_z)\geq0$.

\

If $-D.(L_x+L_y+L_z)=D.e_5$, by the above computation one has (here
for simplicity, we let $d=D.l,m_i=D.e_i,m_{ij}=D.e_{ij}$):
\begin{equation*}
2m_i-\sum_{j\neq i}m_{ij}=0, \quad d-m_5-m_i+m_{i5}=0,\quad
m_{jk}+m_{kl}+m_{jl}+m_{i5}=d.
\end{equation*}
It follows that:
\begin{equation}\label{rel1}
\quad m_{ij}+m_{ik}+m_{il}=d-m_5+m_i,
\end{equation}
\begin{equation}\label{rel2}
m_{jk}+m_{kl}+m_{jl}=2d-m_5-m_i.
\end{equation}
Adding up all relations (\ref{rel1}) and (\ref{rel2}), one has:
$$2\sum_{i,j=1,\ldots 4}m_{ij}=4d-4m_5+\sum_{i=1}^4m_i,\quad
2\sum_{i,j=1,\ldots 4}m_{ij}=8d-4m_5-\sum_{i=1}^4m_i.$$

It follows that $\sum_{i=1}^4m_i=2d$. But by assumption $m_i\geq
m_5$ for all $i$, hence $m_5\leq d/2$. As $0\leq m_{i5}=m_i+m_5-d$
it follows that $m_i\geq d-m_5\geq d/2$. Since $\sum_{i=1}^4m_i=2d$
it follows that $m_i=d/2$, $m_{i5}=0$. Moreover,
$m_{ij}+m_{ik}+m_{il}=d$. As $d>0$ it follows that $m_{ij}>0$ for
some $i,j\in\{1,\ldots,4\}$. This contradicts the assumption in the
Main Claim. \ep

\section{Multiplicity estimates}\label{multiplicities}

Let $l$ be the unique line in $\PP^3$ that passes through $p_5$ and
intersecting lines $l_{13}$ and $l_{24}$ (the other cases are
similar). Let $L$ be the proper transform of $l$ in $\M$.

\begin{prop}\label{multiplicity estimates}
Let $D=dH-\sum m_i E_i-m_{ij}E_{ij}$ be an effective divisor on
$\M$. Let $m$ be the multiplicity of $D$ along $L$. Then
$$m\geq m_5+m_{13}+m_{24}-d.$$
\end{prop}

\bp Let $\rho:X\ra\M$ be the blow-up of $\M$ along $L$ and let $E$
be the exceptional divisor. Let $\tilde{D}$ be the proper transform
of $D$. Then $\rho^*D=\tilde{D}+mE$. Restricting to $E$, one has:
\begin{equation}\label{m}
(\rho^*D)_{|E}=\tilde{D}_{|E}+mE_{|E}.
\end{equation}

Let $N$ be the normal bundle of $L$ in $\M$. Let $N_{l|\PP^3}$ be
the normal bundle of $l$ in $\PP^3$. If $l'$ is the proper transform
of $l$ in the blow-up $X$ of $\PP^3$ along $p_1,\ldots,p_5$, let
$N'$ be the normal bundle of $l'$ in $\M$. One has:
\begin{equation}\label{normal bundle}
N_{l|\PP^3}=\O(-1)\oplus\O(-1),\quad
N'=\pi^*N_{l|\PP^3}(-E_5)=\O\oplus\O.
\end{equation}
It is easy to see that $\deg(N)=\deg(N')-2=-2$. In fact we have the
following:
\begin{claim}
$N=\O(-1)\oplus\O(-1)$.
\end{claim}

\bp Note that one could obtain $\M$ by blowing up $\PP^3$ first
along the points $p_1,\ldots,p_4$, then the proper transforms of the
lines $l_{13}$ and $l_{24}$, then the point $p_5$ and the proper
transforms of the lines $l_{ij}$, for all $ij\neq13, 24$. Let $\La$
be the plane in $\PP^2$ spanned by the line $l$ and $l_{13}$. Then
the proper transform $\tilde{\La}$ of $\La$ in $\M$ is the blow-up
of $\La\cong\PP^2$ along $p_1,p_3,p_5,q$, where $q=l_{24}\cap\La$.
If $N_{L|\tilde{\La}}$ is the normal bundle of $L$ in $\tilde{\La}$
and $N_{\tilde{\La}|\M}$ is the normal bundle of $\tilde{\La}$ in
$\M$, one has an exact sequence:
\begin{equation}\label{exseq}
0\ra N_{L|\tilde{\La}}\ra N\ra (N_{\tilde{\La}|\M})_{|L}\ra0.
\end{equation}
It is easy to see that $ N_{L|\tilde{\La}}=\O(-1)$. Since $\deg(N)=-2$
and $\O(-1)$ is a subbundle of $N$ (the quotient is a line bundle), it
follows that $N=\O(-1)\oplus\O(-1)$.
\ep

Then $E=\PP(N)\cong\PP^1\times\PP^1$. Let $p:E\ra\l=\PP^1$ be the
restriction of $\rho$ to $E$. Let $q:\PP^1\times\PP^1\ra\PP^1$ be
the other projection. Then $$E_{|E}=\O_E(-1)=q^*\O(-1)\otimes
p^*\O(-1).$$

Note that $(\rho^*D)_{|E}=p^*(D_{|L})$ and $D_{|L}=\O(a)$, where we
let $a=D.L$. One has:
$$H.L=E_5.L=E_{13}.L=E_{24}.L=1, E_i.L=E_{ij}.L=0,
\hbox{ for all other indices }i,j.$$ It follows that
$a=d-m_5-m_{13}-m_{24}$. From (\ref{m}) one has:
$$\tilde{D}_{|E}=p^*\O(a+m)\otimes q^*\O(m).$$
Since $\tilde{D}_{|E}$ is effective, it follows that $m\geq
-a=m_5+m_{13}+m_{24}-d$. \ep

\section{Proof of Lemma \ref{Cox of Y}}\label{proof Cox of Y}

Recall that $Y$ is the blow-up of $\PP^2$ along
$q_1,\ldots,q_4,x,y,z$.
\begin{figure}[h]
\centerline{ \psfig{figure=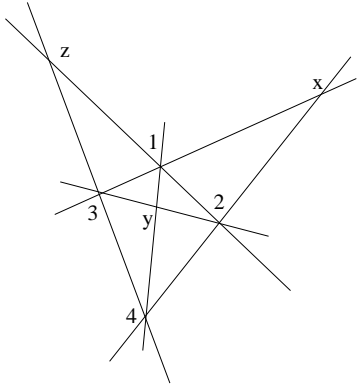,height=2in} }
\smallskip
\centerline{
\parbox{4.5in}{\caption[]{{\small The configuration of the points
$q_1,q_2,q_3,q_4,x,y,z$ }}} }
\end{figure}

Let
$$D=d\H-\sum_{i=1}^4 m_i\E_i-m_x\E_x-m_y\E_y-m_z\E_z,$$ be a
divisor on $Y$. Assume $D$ is effective and let $s$ be
a section in $\HH^0(Y,D)$.

We show that $s$ is generated by distinguished sections on $Y$ by
induction on $d$. Let $\l'_{ij}$ (respectively
$\l_{xy}$,$\l_{yz}$,$\l_{xz}$) be the proper transforms in $Y$ of
the lines $\l_{ij}$ (respectively $\overline{xy}$,
$\overline{yz}$,$\overline{xz}$). We may assume $D.C\geq0$ for $C$
among the classes:
$$\l'_{ij},\l_{xy},\l_{xz},\l_{yz},\E_i,\E_x,\E_y,\E_z.$$

This is because if $D.C<0$ then $s=x_Cs'$, where $s'\in\HH^0(Y,D-C)$
and $x_C$ is a generator of $\HH^0(Y,C)$, and $s'$ is generated by
distinguished sections by induction. Hence, we assume:
\begin{gather*}\tag{$*$}\label{@}
d\geq m_i+m_j+m_{\chi(ij)},\quad
d\geq m_x+m_y,\quad d\geq m_x+m_z,\quad d\geq m_y+m_z,\\
d\geq m_i\geq0,\quad d\geq m_x\geq0,\quad d\geq m_y\geq0,\quad d\geq
m_z\geq0.
\end{gather*}

If $d=0$ then it follows by (\ref{@}) $D=0$. Assume $d>0$. We may
assume without loss of generality that 
% correction 3
$$m_4\leq m_1\leq m_2\leq m_3.$$
% end correction 3

Consider the restriction map:
\begin{equation*}
r:\HH^0(Y,D)\ra\HH^0(\E_4,D_{|\E_4})=\HH^0(\PP^1,\O(m_4)).
\end{equation*}

It is enough to show that we may lift any $t\in\HH^0(\PP^1,\O(m_4))$
to a section in $\HH^0(Y,D)$ generated by distinguished sections on
$Y$. This is because by the same argument as in Section \ref{plan},
if $s,s'$ are sections in $\HH^0(Y,D)$ are such that $r(s)=r(s')$, then
$s-s'$ is in $\HH^0(Y,D-\E_4)$ and we are done by induction.

% correction 4
Let $t_i$ be the restriction in $\HH^0(\PP^1,\O(1))$ of the section
$s_{i4}$ corresponding to $\l'_{i4}$. Any two of $t_1,t_2,t_3$
generate $\HH^0(\PP^1,\O(1))$. In particular, it is enough  to lift
$t={t_1}^k{t_3}^{m_4-k}$ ( for any $0\leq k\leq m_4$) to a combination of
distinguished sections. We lift $t_i$ to $s_{i4}$, hence $t$ to
$s_{14}^k s_{34}^{m_4-k}$ (a section of
$D'=k\l'_{14}+(m_4-k)\l'_{34}$). Let:
\begin{gather*}
\De=D-D'=(d-m_4)\H-(m_1-k)\E_1-m_2\E_2-(m_3-m_4+k)\E_3-\\
-m_x\E_x-(m_y-k)\E_y-(m_z-m_4+k)\E_z.
\end{gather*}

\begin{claim}\label{effective 1}
There is a section $u\in\HH^0(\De)$, generated by distinguished sections and such that 
$u_{|\E_4}\in\HH^0(\E_4,\O)$ is non-zero.
\end{claim}

Assuming Claim \ref{effective 1}, we lift $t$ to $us_{14}^k
s_{24}^{m_4-k}$ and we are done. 

\bp[Proof of Claim \ref{effective 1}]  
Let 
\begin{gather*}
\De'=\De-(m_1-k)\l'_{12}=(d-m_4-m_1+k)\H-(m_2-m_1+k)\E_2-\\
-(m_3-m_4+k)\E_3-m_x\E_x-(m_y-k)\E_y-(m_z-m_1-m_4+2k)\E_z.
\end{gather*}
Note that since $k\leq m_4\leq m_1$ and since a section corresponding to $\l'_{12}$ has non-zero restriction
to $\E_4$, it is enough to show that there is a section $u'\in\HH^0(\De')$, generated by distinguished sections 
and such that  $u'_{|\E_4}\in\HH^0(\E_4,\O)$ is non-zero. 

\noindent\underline{\bf Case when $m_y-k<0$.}
Let 
\begin{gather*}
\De''=\De+(m_y-k)\E_y=(d-m_4-m_1+k)\H-(m_2-m_1+k)\E_2-\\
-(m_3-m_4+k)\E_3-m_x\E_x-(m_z-m_1-m_4+2k)\E_z.
\end{gather*}
It is enough to show that there is a section $u''\in\HH^0(\De'')$, generated by distinguished sections and such that 
$u''_{|\E_4}\in\HH^0(\E_4,\O)$ is non-zero. Since $\De''$ is a divisor on the blow-up of $\PP^2$ along 
the points $q_2, q_3, x, z$, it follows from Lemma \ref{eff cone} (a direct check shows that all inequalities 
(\ref{*}) hold; use $k\leq m_4\leq m_i$ and (\ref{@})) and Lemma \ref{non-vanishing} applied to the lines $\ov{q_2,x}$ 
and $\ov{q_3,z}$, that there is a section $u''\in\HH^0(\De'')$, generated by distinguished sections and not containing $q_4$ in its zero-locus.

\

\noindent\underline{\bf Case when $m_y-k\geq 0$.} 
Denote 
$$N_1=m_1+m_4+m_x+m_y-d-2k,$$
$$N_2=2d-m_2-m_3-m_x-m_z-2k.$$
\begin{claim}\label{valid}
$N_1\leq N_2$, $0\leq N_2$, $N_1\leq m_y-k$.
\end{claim}

\bp[Proof of Claim~\ref{valid}]
We have
\begin{equation*}
N_2-N_1=(d-m_1-m_2-m_z)+(d-m_3-m_4-m_x)+(d-m_x-m_y)\geq0,
\end{equation*}
using (\ref{@}) and $m_4\leq m_i$. Similarly, as $0\leq k\leq m_4$, we have
$N_2\geq0$ and $N_1\leq m_y-k$ (using (\ref{@}) and $m_4\leq m_i$).
\ep

By Claim \ref{valid}, we may choose $\al,\be\geq0$ be integers such that
$\al+\be=m_y-k$ and $N_1\leq\al\leq N_2$. Let 
\begin{gather*}
\De''=\De-\al\l_{xy}-\be\l'_{23}=(d-m_1-m_4-m_y+2k)\H-(m_2-m_1+k-\be)\E_2-\\
-(m_3-m_4+k-\be)\E_3-(m_x-\al)\E_x-(m_z-m_1-m_4+2k)\E_z.
\end{gather*}
Since $\l_{xy}$ and $\l'_{23}$ have non-zero restriction to $\E_4$, it is enough to find 
$u\in\HH^0(\De'')$ such that $u_{|\E_4}\neq0$. As before, since $\De''$ is a divisor on the blow-up of $\PP^2$ along 
the points $q_2, q_3, x, z$, it follows from Lemma \ref{eff cone}  and  Lemma \ref{non-vanishing} applied to the lines 
$\ov{q_2,x}$  and $\ov{q_3,z}$, that there is a section $u''\in\HH^0(\De'')$, 
generated by distinguished sections and not containing $q_4$ in its zero-locus. 
All inequalities follow in a straightforward way from (\ref{@}) and $m_4\leq m_i$, except for:
\bi
\item $\De''.(\H-\E_x)\geq0$, (equivalent to $\al\geq N_1$)
\item  $\De''.(\H-\E_3)\geq0$ (use that $m_1\leq m_2$) 
\item $\De''.(2\H-\E_2-\E_3-\E_x-\E_z)\geq0$ (equivalent to $\al\leq N_2$)
\ei
\ep
% end correction 4

\begin{lemma}\label{eff cone}
Let $Z$ be the blow-up of $\PP^2$ along points $q_1,\ldots,q_4$ (no
three collinear). One has $\HH^0(D)\neq0$ for a divisor
$D=d\H-\sum_{i=1}^4m_i\E_i$ if and only if:
\begin{equation}\label{*}
d\geq 0,\quad d-m_i\geq0\quad 2d-\sum_{i=1}^4m_i\geq0.
\end{equation}
The Cox ring $\Cox(Z)$ is generated by sections corresponding to the
lines $\l_{ij}$ and the exceptional divisors $\E_i$.
\end{lemma}

\bp It is a well known result that the Cox ring $\Cox(Z)$ is
generated by sections corresponding to the lines $\l_{ij}$ and the
exceptional divisors $\E_i$, see for example \cite{BP}.

If $D$ is an effective divisor, then clearly, the inequalities
(\ref{*}) hold. Conversely, assume (\ref{*}) hold. We write $D$ as
an effective combination of the classes of the lines
$\l_{ij}=\H-\E_i-\E_j$ and the exceptional divisors $\E_i$. Consider
the table with $2$ rows and $d$ columns filled with $\E_i$'s in the
following way. Start in the upper left corner and write $m_1$
$\E_1$'s in the first row. Then write $m_2$ $\E_2$'s passing to the
second row if necessary, and so on. Fill the remaining entries with
zeros. For example, if $D=5\H-3\E_1-3\E_2-2\E_3-\E_4$:
$$\begin{matrix}
\E_1&\E_1&\E_1&\E_2&\E_2,\cr \E_2&\E_3&\E_3&\E_4&0.\cr
\end{matrix}
$$
Our conditions guarantee that all entries of a given column are
different. Therefore $D$ is the sum of classes $H-(\E_i+\E_j)$, one
for each column, where $\E_i,\E_j$ are the entries of the column. In
the example above:
$$D=(\H-\E_1-\E_2)+(\H-\E_1-\E_3)+(\H-\E_1-\E_3)+(\H-\E_2-\E_4)+(\H-\E_2).$$\ep

% added (correction 5)
\begin{lemma}\label{non-vanishing}
In the notations of Lemma \ref{eff cone}, let $D$ be a divisor such that $\HH^0(D)\neq0$.
Let $q$ be the intersection point of the lines $\ov{q_1q_2}$ and $\ov{q_3q_4}$. The linear system $|D|$
does not contain $q$ as a base point if and only if 
$$D.(\H-\E_1-\E_2)\geq0,\quad D.(\H-\E_3-\E_4)\geq0.$$
\end{lemma}

\bp
The conditions are clearly necessary. It is enough to show that 
$D$ can be written as an effective combination of lines 
$\l_{ij}$ ($\l_{ij}\neq\l_{12},\l_{34}$) and the exceptional divisors $\E_i$.
Let $$D=\sum k_{ij}\l_{ij}+\sum k_i\E_i, \quad k_{ij}, k_i\geq0.$$

Assume $k_{12}>0$. Note that the only generators $E$ of $\Cox(Z)$ with the property that 
$E.\l_{12}>0$ are $\l_{34}, \E_1, \E_2$.  
Since $D.\l_{12}\geq0$, it follows that one of $k_{34}, k_1, k_2>0$.
If $k_1>0$ we may replace $\l_{12}+\E_1$ with a divisor in the pencil $|\H-\E_2|$ that does not contain $\l_{12}$
(for example $\l_{23}+\E_3$). The case $k_2>0$ is similar. If $k_{34}>0$, we replace $\l_{12}+\l_{34}$ with, for example,
$\l_{13}+\l_{24}$. The case when $k_{34}>0$ is similar. At the end of this process, we have $k_{12}=k_{34}=0$.
\ep
% end correction 5

\section{Inequalities for the effective cone of $X$}\label{LM}

Let $X$ be the iterated blow-up of $\PP^3$ in points
$p_1,\ldots,p_4$ (in linearly general position) and proper
transforms of lines $l_{ij}$ ($i,j=1,\ldots 4$, $i\neq j$). Let
$E_i,E_{ij}$ be the exceptional divisors. Let $l$ be the class on
$X$ of the proper transform of a general line in $\PP^3$. Let $e_i$
be the class of (the proper transform of) a general line in $E_i$.
Let $e_{ij}$ be the class of a fiber of the $\PP^1$-bundle
$E_{ij}\ra l_{ij}$.

\begin{notn}
For $\{i,j,k,l\}=\{1,2,3,4\}$ let:
$$C_{ij}=2l-e_{ij}-e_k-e_l,$$
$$C_{i;j}=2l-e_{ik}-e_{il}-e_j,$$
$$C_i=2l-e_{ij}-e_{ik}-e_{il},$$
$$B=3l-\sum_{i=1}^4e_i,$$
$$B_i=3l-2e_i-e_{jk}-e_{jl}-e_{kl}.$$
\end{notn}

% correction 6 
\begin{lemma}\label{LMspaces}
Let $D$ be a divisor on $X$. Then $D$ is an effective sum (with integer, non-negative coefficients) 
of boundary divisors $\La_{ijk},E_{ij},E_i$ (in particular, $\HH^0(D)\neq0$)
if $D.C\geq0$ for all $C$ in the list below (for all
$\{i,j,k,l\}=\{1,2,3,4\}$): $(1)$ $l$; $(2)$ $l-e_i$; $(3)$
$l-e_{ij}$; $(4)$ $l-e_{ij}-e_{kl}$; $(5)$ $C_{ij}$; $(6)$
$C_{i;j}$; $(7)$ $C_i$; $(8)$ $B$; $(9)$ $B_i$, and moreover, if one
has $D.(l-e_{ij}-e_{kl})>0$ for some $i,j,k,l$.
\end{lemma}

It is easy to see that each of the classes $C$ in $(1)$-$(9)$ in 
Lemma \ref{LMspaces} cover a dense set of $X$; hence, for any 
effective divisor $D$ one has $D.C\geq0$, i.e., $C$ is a \emph{nef} curve.
% end correction 6

\begin{rmk}
It is a standard fact that the divisor $D$ is in the convex hull of
the effective divisors $\La_{ijk},E_{ij},E_i$ (where $\La_{ijk}$ is
the proper transform of the plane $\overline{p_ip_jp_k}$) if and
only if inequalities $(1)$-$(9)$ hold. 
However, the extra condition
of having at least one strict inequality in $(4)$ is necessary for
$\HH^0(D)\neq0$, as the following example shows: if
$D=2H-\sum_{i\neq j}E_{ij}$ then it is easy to see that $\HH^0(D)=0$
($\HH^0(2D)\neq0$) and $D$ satisfies all of $(1)$-$(9)$.
\end{rmk}

\begin{obs}\label{negatives}
If $$D=dH-\sum m_iE_i-\sum m_{ij}E_{ij},$$ is such that $d\geq0$,
$d\geq m_i,m_{ij}$ (for all $i,j$) and there is an $i$ such that
$m_i\leq0$ and $m_{ij}\leq0$ for all $j\neq i$, then
% correction 7 
$D$ is an effective sum of boundary, 
% end correction 7
as one has:
$$D=d\La_{jkl}+\sum_{j\neq i}(d-m_j)E_j+
\sum_{u,v\neq i}(d-m_{uv})E_{uv}+(-m_i)E_i+\sum_{j\neq
i}(-m_{ij})E_{ij}.$$
\end{obs}

\bp[Proof of Lemma~\ref{LMspaces}] Let $D=dH-\sum m_iE_i-\sum
m_{ij}E_{ij}$. One has: $$d=D.l,\quad m_i=D.e_i,\quad
m_{ij}=D.e_{ij}.$$

We do an induction on $d$. If $d=0$ then from $(2)$ and $(3)$
$m_i,m_{ij}\leq0$
% correction 8
and we are done by Observation \ref{negatives}.
% end correction 8
Assume $d>0$. We show
that there are $i,j,k$ such that $D'=D-\La_{ijk}$ also satisfies
$(1)$-$(9)$ and hence, we are done by induction.

Note that $D'=D-\La_{ijk}$ (for any $i,j,k$) always satisfies
$(1)$,$(4)$,$(5)$,$(8)$,$(9)$. Moreover, one has at least one strict
inequality in $(4)$. Inequality $(2)$ fails for $D'$ if and only if
$m_l=d$, where $l\neq i,j,k$, and $(3)$ fails for $D'$ if and only
if one of $m_{il},m_{jl},m_{kl}$ equal $d$. Inequality $(6)$ fails
for $D'$ if and only if one has:
$$2d=m_{li}+m_{lj}+m_k,$$
(or the similar inequalities for a permutation of indices $i,j,k$).
Inequality $(7)$ fails for $D'$ if and only if
$m_{il}+m_{jl}+m_{kl}\in\{2d-1,2d\}$.

\

\noindent\underline{{\bf Case I}: $m_{ij}=d$ for some $i,j$}. We may
assume $d=m_{12}$. From $(4)$, $m_{34}\leq0$.

\

\noindent\underline{Case 1: $m_i=d$ for $i\in\{3,4\}$}. We may
assume $m_4=d$. Then by $(5)$ one has that $m_3\leq0$ and by $(6)$
one has $m_{13},m_{23}\leq0$ and we are done by Observation
\ref{negatives}.

\

\noindent\underline{Case 2: Assume $m_3<d,m_4<d$.}

We may assume that $m_{13}$ is the largest among
$m_{13},m_{14},m_{23},m_{24}$.
\begin{claim}\label{strict}
One has $m_{14},m_{24}<d$.
\end{claim}
\bp Assume $m_{i4}=d$ for $i=\{1,2\}$. Since by assumption
$m_{i4}\leq m_{13}$ and $m_{13}\leq d$ one has $m_{13}=d$. If
$m_{14}=d$, since $m_{12}=d$, one has a contradiction with $(7)$. If
$m_{24}=d$ one contradicts $(4)$. \ep

\begin{claim}\label{D minus plane}
The divisor $D'=D-\La_{123}$ satisfies $(1)$-$(9)$.
\end{claim}

\bp Inequality $(2)$ holds, as  $m_4<d$. Since by Claim \ref{strict}
$m_{14},m_{24}<d$ and since $m_{34}\leq0$, inequality $(3)$ holds.
If $(7)$ is not satisfied, i.e., $m_{14}+m_{24}+m_{34}=2d-1$ or
$2d$, one has a contradiction, as $m_{34}\leq0$ and (by Claim
\ref{strict}) $m_{14},m_{24}<d$. If $(6)$ is not satisfied for $D'$
then one has:
\begin{equation}\label{eq1}
2d=m_{i4}+m_{j4}+m_k,
\end{equation}
for some $\{i,j,k\}=\{1,2,3\}$. If $k=3$: by $(6)$ one has $2d\geq
m_{12}+m_{24}+m_3$. Since $m_{12}=d$ one has $d\geq m_{24}+m_3$ and
hence, by (\ref{eq1}), $m_{14}=d$. This contradicts Claim
\ref{strict}. If $k\in\{1,2\}$ (say $i=3$): one has
$2d=m_{34}+m_{j4}+m_k$. Since $m_{34}\leq0$, $m_{j4}<d$ (Claim
\ref{strict}) and $m_k\leq d$, this is a contradiction. \ep

\noindent\underline{{\bf Case II}: $m_i=d$ for some $i$, $m_{ij}<d$
for all $i,j$}. We may assume that $d=m_4$. We may also assume that
$d\geq m_3\geq m_2\geq m_1$. By $(9)$ and $(5)$ one has:
\begin{equation}\label{*1}
m_{12}+m_{13}+m_{23}\leq d,
\end{equation}
\begin{equation}\label{*2}
m_{ij}+m_k\leq d,
\end{equation}
where $\{i,j,k\}=\{1,2,3\}$.

\noindent\underline{Case 1: $m_{23}>0$.} By (\ref{*1}) one has:
\begin{equation}\label{*4}
m_{12}+m_{13}<d.
\end{equation}

\begin{claim}
The divisor $D'=D-\La_{234}$ satisfies $(1)$-$(9)$.
\end{claim}

\bp Inequality $(3)$ is satisfied by the assumption $m_{ij}<d$.
Inequality $(2)$ is not satisfied if and only if one has $m_1=d$. If
$m_1=d$, it follows from the assumptions that $m_2=m_3=d$. As
$m_4=d$, one has a contradiction with $(8)$. If $(7)$ is not
satisfied, i.e., $m_{12}+m_{13}+m_{14}=\{2d-1,2d\}$, one has a
contradiction with $m_{14}<d$ and $m_{12}+m_{13}<d$ (\ref{*4}). If
$(6)$ is not satisfied then one has:
\begin{equation}\label{eq2}
2d=m_{1i}+m_{1j}+m_k,
\end{equation}
for some $\{i,j,k\}=\{2,3,4\}$. If $k=4$: since $m_4=d$ one has from
(\ref{eq2}) $d=m_{12}+m_{13}$ which contradicts (\ref{*4}).  If
$k\in\{2,3\}$ (say $i=4$) one has $2d=m_{14}+m_{1j}+m_k$ for
$\{k,j\}=\{2,3\}$. But $m_{14}<d$ and $m_{j1}+m_k\leq d$ (\ref{*2}).
This is a contradiction. \ep

\noindent\underline{Case 2: $m_{23}\leq0$.} If $m_2=d$ then it
follows from the assumptions that $m_3=d$. As $m_4=d$, one has from
$(8)$ that $m_1\leq0$. It follows from $(5)$ that $m_{1i}\leq0$ for
all $i=2,3,4$. Then we are done by Observation \ref{negatives}.
Hence, we may assume $m_2<d$.

\begin{claim}
The divisor $D'=D-\La_{134}$ satisfies $(1)$-$(9)$.
\end{claim}

\bp Inequality $(2)$ is satisfied as $m_2<d$. Inequality $(3)$ is
satisfied by assumption. If $(7)$ is not satisfied, i.e.,
$m_{12}+m_{23}+m_{24}\in\{2d-1,2d\}$, one has a contradiction with
$m_{24}<d$ and $m_{12}+m_{23}\leq m_{12}<d$. If $(6)$ is not
satisfied then:
\begin{equation}\label{eq3}
2d=m_{2i}+m_{2j}+m_k,
\end{equation}
for some $\{i,j,k\}=\{1,3,4\}$. If $k=4$: since $m_4=d$ one has from
(\ref{eq3}) $d=m_{12}+m_{23}$. But $m_{12}<d$ and $m_{23}\leq0$.
This is a contradiction. If $k\in\{1,3\}$ (say $i=4$) one has
$2d=m_{24}+m_{j2}+m_k$ for $\{k,j\}=\{1,3\}$. But $m_{24}<d$ and
$m_{j2}+m_k\leq d$ (\ref{*2}). This is a contradiction. \ep

\noindent\underline{{\bf Case III}: $m_i<d$, $m_{ij}<d$ for all
$i,j$.}

\

By Claim \ref{l;i} we may assume $D.C_{i;j}>0$ for $i=1,2,3$, and
all $j\neq i$.

\begin{claim}
One of $D_1=D-\La_{234}$, $D_2=D-\La_{134}$, $D_3=D-\La_{124}$
satisfies all the inequalities $(1)-(9)$.
\end{claim}

\bp  Inequalities $(2)$,$(3)$ follow from the assumptions. If $(6)$
is not satisfied for $D_i$ then $D.C_{i;j}$, for some $j\neq i$,
which we assume does not happen. Hence, $(6)$ is satisfied for all
$D_i$. If $(7)$ fails for all $D_i$, then one has for all
$i\in\{1,2,3\}$:
\begin{equation}\label{eq4}
\quad m_{ij}+m_{ik}+m_{il}\geq 2d-1.
\end{equation}

Adding up (\ref{eq4}) for  $i=1$ and $i=2$ one has:
\begin{equation}\label{eq5}
2m_{12}+(m_{13}+m_{24})+(m_{14}+m_{23})\geq 4d-2.
\end{equation}

From $(4)$ $d\geq m_{ij}+m_{kl}$. As $m_{12}<d$, it follows from
(\ref{eq5}) that $m_{13}+m_{24}=m_{14}+m_{23}=d$. Similarly, by
adding (\ref{eq5}) for $i=1$ and $i=3$ one has $m_{12}+m_{34}=d$.
This contradicts our assumption that one of the inequalities in
$(4)$ is strict. \ep

\begin{claim}\label{l;i}
There are at least three indices $i\in\{1,2,3,4\}$ such that
$D.C_{i;j}>0$ for all $j\neq i$.
\end{claim}

\bp Assume $D.C_{i;j}=0$, for some $i,j$. We may assume without loss
of generality that $D.C_{1;2}=0$.
\begin{equation*}
2d=m_{13}+m_{14}+m_2.
\end{equation*}

We claim that for all $i\in\{2,3,4\}$ one has $D.C_{i;j}>0$ for all
$j\neq i$. This follows from $D.(C_{1;2}+C_{i;j})>0$ for all
$i\in\{2,3,4\}$, $j\neq i$. This is because:
$$D.(C_{1;2}+C_{3;j})=D.C_{3k}+D.(l-e_{23}-e_{14})+D.e_{13} \quad (\{j,k\}=\{1,4\}).$$
It follows from $(5)$, $(4)$ and $m_{13}<d$ that
$D.(C_{1;2}+C_{3;j})>0$. Similarly:
$$D.(C_{1;2}+C_{3;2})=D.B_2+D.e_{13}.$$
By $(9)$ and $m_{13}<d$, $D.(C_{1;2}+C_{3;2})>0$. By symmetry,
$D.C_{4;j}>0$, for all $j\neq4$.

If $\{j,k\}=\{3,4\}$, one has:
$$D.(C_{1;2}+C_{2;j})=D.C_{1k}+D.(l-e_{1j}-e_{2k})+D.e_{12}.$$
From $(5)$, $(4)$ and $m_{12}<d$ one has $D.(C_{1;2}+C_{2;j})>0$.
Similarly:
$$D.(C_{1;2}+C_{2;1})=D.(2l-e_1-e_2)+D.(l-e_{13}-e_{24})+D.(l-e_{14}-e_{23}).$$
From $(4)$ and $m_1,m_2<d$ one has $D.(C_{1;2}+C_{2;1})>0$.

\ep \ep

\section{Restrictions to generators}\label{restr to generators}

Let $\pi':\M'\ra\PP^3$ be the blow-up along $p_1,\ldots,p_5$ and let
$E'_i$ be the corresponding exceptional divisors. Let $\pi:\M\ra\M'$
be the blow-up of the proper transforms of the lines $l_{ij}$. In
what follows, we compute the classes of the restrictions of an
arbitrary divisor $D$ on $\M$ to the divisors
$E_i,E_{ij},\La_{ijk},Q_{(ij)(kl)}$.

\subsection{Restrictions to $E_i$.}\label{section restr to E_i}
The divisor $E_i$ is the inverse image $\pi^{-1}(E'_i)$. By Fact
\ref{blowups} the divisor $E_i$ is the blow-up of $E'_i\cong\PP^2$
along the $4$ points corresponding to the directions of the lines
$l_{ij}$, for $j\neq i$. Denote by $\E_j$ the corresponding
exceptional divisors. Denote by $\H$ the hyperplane class on $E_i$.
One may easily see the following:
\begin{equation}\label{basic restr to E_i}
H_{|E_i}=0,\quad {E_i}_{|E_i}=-\H,\quad {E_j}_{|E_i}=0\quad
{E_{ij}}_{|E_i}=\E_j, \quad {E_{jk}}_{|E_i}=0,
\end{equation}
where $j, k\neq i$, $j\neq k$. This is clear from Fact
\ref{blowups}.

\begin{fact}\cite[Prop. IV.21,p.167]{EH}\label{blowups}
Let $Y$ and $Z$ be closed subschemes
in a scheme $X$ and let $\tilde{X}$ be the blow-up of $X$ along $Z$.
Let $E$ be the exceptional divisor. The proper transform $\tilde{Y}$
of $Y$ is the blow-up of $Y$ along the scheme theoretic intersection
$Y\cap Z$ and the exceptional divisor is $\tilde{Y}\cap E$. In
particular, if $Z$ is contained in $Y$, the scheme $\tilde{Y}$ is
the blow-up of $Y$ along $Z$.
\end{fact}

Consider an arbitrary divisor $D$ on $\M$:
\begin{equation}\label{D}
D=dH-\sum_i m_iE_i-\sum_{i,j}m_{ij}E_{ij}, \hbox{ where } d, m_i,
m_{ij}\in\ZZ.
\end{equation}

It follows from (\ref{basic restr to E_i}) that the restriction of
$D$ to $E_i$ is given by:
\begin{equation}\label{restr to E_i}
D_{|E_i}=m_i\H-\sum_{j\neq i}m_{ij}\E_j.
\end{equation}

\begin{lemma}
The divisor $D_{|E_i}$ is an effective divisor if and only if
\begin{equation}\label{ineq E_i}
m_i\geq 0,\quad m_i\geq m_{ij} \quad (j\neq i),\quad 2m_i\geq
\sum_{j\neq i}m_{ij}.
\end{equation}
\end{lemma}

\bp This is Lemma \ref{eff cone} applied to (\ref{restr to E_i}).\ep

\subsection{Restrictions to $E_{ij}$.} The normal bundle $N$ of the
proper transform of the line $l_{ij}$ in $\M'$ is given by:
$$N=(\pi^*N_{l_{ij}|\PP^3})(-E_i-E_j)=\O(-1)\oplus\O(-1).$$

The divisor $E_{ij}=\PP(N)$ is isomorphic to
$\PP(\O\oplus\O)=\PP^1\times\PP^1$. Let $p_1:\PP^1\times\PP^1\ra
\PP^1$ be the projection map given by the blow-up map $E_{ij}\ra
l_{ij}=\PP^1$ and let $p_2$ be the other projection. Since
$\O(E_{ij})_{|E_{ij}}=\O_{\PP(N)|\PP^1}(-1)$ and
$$\O_{\PP(N)|\PP^1}(-1)\cong\O_{\PP(\O\oplus\O)|\PP^1}(-1)\otimes
p_1^*\O(-1),$$ it follows that:
\begin{equation*}
{E_{ij}}_{|E_{ij}}=p_1^*\O(-1)\otimes p_2^*\O(-1).
\end{equation*}
Moreover, one may easily see, for all distinct $i,j,k,l$:
\begin{equation}\label{basic restr to E_ij}
H_{|E_{ij}}={E_i}_{|E_{ij}}=p_1^*\O(1), \quad
{E_k}_{|E_{ij}}=0,\quad {E_{kl}}_{|E_{ij}}={E_{ik}}_{|E_{ij}}=0.
\end{equation}
It follows from (\ref{basic restr to E_ij}) that the restriction of
$D$ in (\ref{D}) to $E_{ij}$ is given by:
\begin{equation*}
D_{|E_{ij}}=p_1^*\O(d-m_i-m_j+m_{ij})\otimes p_2^*\O(m_{ij}).
\end{equation*}
Clearly, the divisor $D_{|E_{ij}}$ is an effective divisor if and
only if
\begin{equation}\label{ineq E_ij}
m_{ij}\geq 0,\quad d-m_i-m_j+m_{ij}\geq0.
\end{equation}

\subsection{Restrictions to $\La_{ijk}$.}
Take the case of $\La_{123}$ (the other cases are similar). Let
$\La$ be the plane $\overline{p_1p_2p_3}$. Then $\La_{123}$ is the
proper transform $\tilde{\La}$ of $\La$ in $\M'$. Denote by $\La'$
the proper transform of $\La$ in $\M'$. Let $q$ be the point
$l_{45}\cap\La$. Note that by Fact \ref{blowups}, $\La'$ is the
blow-up $\La=\PP^2$ along $p_1,p_2,p_3$ and $\tilde{\La}$ is
isomorphic to the blow-up of $\La'$ in $q$, i.e., $\tilde{\La}$ is
isomorphic to the blow-up of $\PP^2$ along $p_1,p_2,p_3,q$. Let
$\E_1,\E_2,\E_3,\E_q$ be the exceptional divisors and $\H$ the
hyperplane class. One may easily see that:
\begin{equation*}
H_{|\tilde{\La}}=\H,\quad {E_i}_{|\tilde{\La}}=0 \quad(i=4,5),
\quad{E_{ij}}_{|\tilde{\La}}=0\quad (ij\neq 12,13,23,45).
\end{equation*}
Using Fact \ref{blowups}, one has that:
\begin{equation*}
{E_i}_{|\tilde{\La}}=\E_i \quad (i=1,2,3),\quad
{E_{45}}_{|\tilde{\La}}=\E_q,\quad {E_{ij}}_{|\tilde{\La}}
=\H-\E_i-\E_j\quad (ij\in\{12,13,23\}).
\end{equation*}
It follows that the restriction of $D$ in (\ref{D}) to $\La_{123}$
is given by:
\begin{equation*}
D_{|\La_{123}}=(d-m_{12}-m_{13}-m_{23})\H-
\sum_{\{i,j,k\}=\{1,2,3\}}(m_i-m_{ij}-m_{ik})\E_i-m_{45}\E_q.
\end{equation*}

By permuting indices and applying Lemma \ref{eff cone}, one has the
following:
\begin{lemma}\label{eff}
If the divisor $D_{|\La_{ijk}}$ is effective and
$\{i,j,k,u,v\}=\{1,2,3,4,5\}$ then
\begin{equation}\label{ineq planes}
d\geq m_i+m_{jk},\quad d\geq m_{ij}+m_{ik}+m_{jk}+m_{uv}, \quad
2d\geq m_i+m_j+m_j+m_{uv}.
\end{equation}
\end{lemma}

\subsection{Restrictions to the Keel-Vermeire divisors $Q_{(ij)(kl)}$.}
Take the case of $Q_{(12)(34)}$. There is a unique (smooth) quadric
$Q$ in $\PP^3$ that contains the points $p_1,\ldots,p_5$ and the
lines $l_{13},l_{14},l_{23},l_{24}$. Since $Q_{(12)(34)}$ has class:
$$Q_{(12)(34)}=2H-\sum_i E_i-E_{13}-E_{14}-E_{23}-E_{24},$$
it follows that $Q_{(12)(34)}$ is the proper transform $\tilde{Q}$
of $Q$ in $\M$. Denote by $Q'$ the proper transform of $Q$ in $\M'$.
By Fact \ref{blowups} it follows that $Q'$ is the blow-up of
$Q\cong\PP^1\times\PP^1$ along the points $p_1,\ldots,p_5$. Moreover
$\tilde{Q}\cong Q'$.

Let $F_1$, respectively $F_2$, be the class of the lines in the
ruling of $\PP^1\times\PP^1$ that contains $l_{13}$ and $l_{24}$,
respectively $l_{14}$ and $l_{23}$. Let $\E_1,\ldots\E_5$ be the
exceptional divisors on $\tilde{Q}$, considered as a blow-up of
$\PP^1\times\PP^1$ along $p_1,\ldots,p_5$. By Fact \ref{blowups}:
$$H_{|\tilde{Q}}=F_1+F_2,\quad {E_i}_{|\tilde{Q}}=\E_i,$$
$${E_{ij}}_{|\tilde{Q}}=F_1-\E_i-\E_j \quad(ij=13,24),\quad
{E_{ij}}_{|\tilde{Q}}=F_2-\E_i-\E_j \quad(ij=14,23),$$
$${E_{ij}}_{|\tilde{Q}}=0\quad\hbox{ for all other cases.}$$

It follows that restriction $D_{|\tilde{Q}}$ of the divisor $D$ in
(\ref{D}) to $\tilde{Q}$ is given by:
\begin{gather*}
D_{|\tilde{Q}}=(d-m_{13}-m_{24})F_1+(d-m_{14}-m_{23})F_2-(m_1-m_{13}-m_{14})\E_1-\\
-(m_2-m_{23}-m_{24})\E_2-(m_3-m_{13}-m_{23})\E_3-(m_4-m_{14}-m_{24})\E_4-
m_5\E_5.
\end{gather*}

\noindent {\bf Alternative description of $\tilde{Q}$.} Let
$\rho:\PP^3\setminus\{p_5\}\ra\PP^2$ be the projection from $p_5$
and let $q_i=\rho(p_i)$ ($i=1,\ldots,4$). Let $l_1$ (respectively
$l_2$) be the unique line through $p_5$ in the ruling of $F_1$
(respectively $F_2$).

Let $y$ (respectively $x$) be the image $l_1$ (respectively $l_2$).
The blow-up of $Q=\PP^1\times\PP^1$ in $p_5$ is isomorphic to the
blow-up of $\PP^2$ in $x,y$. Hence, $\tilde{Q}$ is isomorphic to the
blow-up of $\PP^2$ along $p_1,\ldots,p_4,x,y$. Denote by
$\E'_1,\ldots\E'_4,\E_x,\E_y$ be the exceptional divisors
corresponding to the points $p_1,\ldots,p_4,x,y$ and let $\H$ be the
hyperplane class. One may immediately see:
$$\H=\rho^*\O(1)=F_1+F_2-\E_5.$$

Note that lines in the ruling $F_2$ (respectively $F_1$) intersect
$l_1$ (respectively $l_2$), therefore their images in $\PP^2$ all
pass through $y$ (respectively $x$). In particular, the lines
$\overline{q_1q_3}$ and $\overline{q_2q_4}$ intersect in $x$, while
the lines $\overline{q_1q_4}$ and $\overline{q_2q_3}$ intersect in
$y$. Moreover, one has:
$$F_1=\H-\E_x,\quad F_2=\H-\E_y.$$
It follows that:
$$\E_5=\H-\E_x-\E_y, \quad \E_i=\E'_i \quad (i=1,\ldots,4).$$
Hence, the restriction $D_{|\tilde{Q}}$ of the divisor $D$ in
(\ref{D}) to $\tilde{Q}$ is given by:
\begin{gather*}
D_{|\tilde{Q}}=(2d-m_5-m_{13}-m_{14}-m_{23}-m_{24})\H-(m_1-m_{13}-m_{14})\E_1-\\
-(m_2-m_{23}-m_{24})\E_2-(m_3-m_{13}-m_{23})\E_3-(m_4-m_{14}-m_{24})\E_4-\\
-(d-m_5-m_{13}-m_{24})\E_x-(d-m_5-m_{14}-m_{23})\E_y.
\end{gather*}

\begin{lemma}
If the divisor $D_{|\tilde{Q}}$  is effective then
\begin{equation}\label{ineq quadrics}
2d\geq m_5+m_{13}+m_{14}+m_{23}+m_{24},\quad 2d\geq
m_1+m_5+m_{23}+m_{24}.
\end{equation}
\end{lemma}

\bp This follows from Lemma \ref{eff cone}. \ep

\end{document}